\documentclass[paper=a4,abstract=true,numbers=noenddot]{scrartcl}

\usepackage[english]{babel}  
\usepackage{hyperref}
\usepackage[intlimits]{amsmath}
\usepackage{amsthm}
\usepackage{amssymb} 
\usepackage{units} 
\usepackage{amsfonts,mathrsfs,bbm,mathabx}
\allowdisplaybreaks[1]
\usepackage[numbers]{natbib}
\usepackage{graphicx,subcaption}
\usepackage{tikz}
\usetikzlibrary{calc}
\usepackage{pgfplots}
\pgfplotsset{every axis/.append style={thick}}
\pgfplotsset{every axis legend/.append
  style={cells={anchor=west},anchor=west}}
\pgfplotscreateplotcyclelist{mylines}{%
  {red,mark=*},
  {blue,mark=square},
  {green,mark=+},
  {black,mark=o},
  {violet,mark=otimes},
  {brown,mark=triangle},
  {cyan,mark=diamond},
  {orange,mark=10-pointed star},
  {magenta,mark=|}  }
\usepgfplotslibrary{fillbetween}
\usepackage[T1]{fontenc}
\usepackage{mathptmx}
\usepackage[scaled=.92]{helvet}
\usepackage{courier}

\usepackage[a4paper,left=2.5cm,right=2.5cm,top=3cm,bottom=2.5cm]{geometry}
\usepackage{authblk}
\usepackage{amsmath,bm}
\usepackage{dsfont}
\usepackage{upgreek}
\usepackage{mathtools}
\usepackage{xcolor}
\usepackage{algorithm}
\usepackage{algpseudocode}
\theoremstyle{remark}

\makeatletter
\def\toprule{
\noalign{\ifnum0=`}\fi
  \hrule \@height 0\p@ \@width 0pt
  \hrule \@height 0.75\p@ 
  \hrule \@height 5pt \@width 0pt
  \futurelet\@tempa\@xhline}
\def\midrule{\noalign{\ifnum0=`}\fi%
  \hrule \@height 3pt \@width 0pt
  \hrule \@height .5pt 
  \hrule \@height 5pt \@width 0pt
  \futurelet \@tempa\@xhline}
\def\botrule{\noalign{\ifnum0=`}\fi
  \hrule \@height 3pt \@width 0pt
  \hrule \@height 0.75\p@ 
  \hrule \@height 3pt \@width 0pt
  \futurelet\@tempa\@xhline}
\makeatother

\title{Comparison of $\mathcal{H}$-matrix- and FMM-based 3D-ACA for
  a time-domain boundary element method}

\author[1]{M. Schanz}
\author[1]{V. Lakshmi Keshava}
\author[2]{H. De Gersem}
\affil[1]{Institute of Applied Mechanics, Graz University of
  Technology, Technikerstraße 4/II, 8010 Graz, Austria,
  m.schanz@tugraz.at, v.lakshmikeshava@tugraz.at}
\affil[2]{Institute for Accelerator Science and Electromagnetic
  Fields, Technical University of Darmstadt, Schloßgartenstr. 8, 64289
  Darmstadt, Germany, herbert.degersem@tu-darmstadt.de}
\date{}                     
\newcommand{\corr}[1]{{\color{black}#1}}

\definecolor{facecol}{rgb}{0.16, 0.4, 1.}
\definecolor{fibercol}{rgb}{0.,0.65, 0.31}

\newcommand{\op}[1]{\mathcal{#1}}
\newcommand{\kl}[1]{\left(#1\right)}
\newcommand{\Dt}{\Delta t}
\newcommand{\Od}[1]{\operatorname{d}#1}

\newcommand{\vek}[1]{\mathbf{#1}}      
\newcommand{\x}[0]{\vek{x}}            
\newcommand{\y}[0]{\vek{y}}            
\newcommand{\mat}[1]{\mathsf{#1}}      
\newcommand{\cvek}{\mat{c}}
\newcommand{\bvek}{\mat{b}}
\newcommand{\bone}{\mathds{1}} 
\newcommand{\A}{\mat{A}}
\newcommand{\bt}{\mat{b}^{\mathsf{T}}}
\newcommand{\matWeightO}[1]{\mat{W}^{\Dt_n}_{#1}}
\newcommand{\matWeight}[2]{\matWeightO{#1}\!\!\kl{#2}}
\newcommand{\ie}{i.e., }

\newcommand{\eg}{e.g., }
\newcommand{\Eg}{E.g., }

\graphicspath{{pics/}{pics_tmp}}
\usepackage{soul}
\sethlcolor{yellow}

\usepackage{enumitem}
\begin{document}
	
\maketitle
\section*{Abstract}
The homogeneous wave equation is solved by a time-domain boundary element method (BEM) using low-order shape functions for spatial, and the generalised convolution quadrature method (gCQ) by Lopez-Fernandez and Sauter for temporal discretisation. The three-dimensional array of BEM matrices according to a set of complex frequencies in Laplace domain is approximated by generalised Adaptive Cross Approximation (3D-ACA). Its rank is increased adaptively until a prescribed accuracy is reached, relying on a pure algebraic error criterion. The data slices for the selected frequency points are further processed by either the standard $\mathcal{H}$-matrices approach with ACA or by a fast multipole method (FMM). This paper compares both approaches with respect to their demands in storage and computing time. Both techniques are illustrated for calculating the sound scattered by an electric machine, for which the proposed algebraic compression techniques make time-domain BEM feasible for the first time.

\textbf{Keywords}: wave equation; boundary element method; generalised
convolution quadrature; multivariate adaptive cross approximation

\section{Introduction}
Wave propagation problems are common in engineering, \eg in sound radiation of machines, in non-destructive testing or when exploring the underground. These problems are formulated with hyperbolic partial differential equations, \eg in acoustics, electromagnetism or elastodynamics. Even for linear hyperbolic problems, the handling of space and time requires expensive discretisation methods, particularly for scattering problems where an unbounded domain has to be considered. Therefore, many wave propagation problems are preferably solved with the boundary element method (BEM). The basis are
boundary integral equations with retarded potentials as kernels, which build the counterpart to the governing hyperbolic
partial differential equation. The mathematical theory goes back to the beginning of last century
for scalar problems like acoustics by Fredholm, and later for vectorial problems in elasticity by Kupradze~\cite{kupradze79}. The mathematical background of time-dependent boundary integral equations is summarised by Costabel~\cite{costabel04} and extensively discussed in the textbook by Sayas~\cite{sayas2016}.

The first numerical realisation of a time-domain boundary element (BE) formulation originates from Mansur \cite{man} in the 80th of last century. Despite its popularity, it suffers from instabilities (see, \eg\cite{peirce97}). A stable space-time formulation has been published by Bamberger and Ha-Duong~\cite{bambduong}, and has been further explored by the group of Aimi~\cite{aimi08a,aimi12a}. These approaches work directly in time domain. Nonetheless, a transformation to Laplace or Fourier domain results in suitable formulations as well, \eg~\cite{cruri}. Formulations in transformed domains offer efficiency from a memory point of view, since essentially only elliptic problems have to be solved, which corresponds to the size of one time step of the time-domain methods. However, to reconstruct the time-dependent solution, numerous frequency-dependent solutions have to be computed, whose solution by iterative solvers is computationally expensive, especially for higher frequencies. Furthermore, the inverse transformation technique necessitates to choose a suitable set of parameters, which often lacks physical motivation. In contrast, time-domain methods require a huge amount of memory since the matrices have to be stored for each time step (up to a cut-off in 3D). Somehow positioned between transformation and time-domain methods are BE formulations based on the convolution quadrature (CQ) method as proposed by Lubich~\cite{lubich88a,lubich88b}. Such a BE formulation is a true time-stepping method utilising the fundamental solutions and properties in Laplace domain. Applications of the CQ to BE methods can be found, \eg in~\cite{schanz01a,schanz97e}. The generalisation of this seminal technique to variable time-step sizes
has been proposed by L\'opez-Fern\'andez and
Sauter~\cite{lopezfernandez13a, lopezfernandez15b} and is called
the generalised convolution quadrature (gCQ) method. Applications can be found in acoustics with absorbing boundary conditions~\cite{sauterschanz17a} and in thermoelasticity~\cite{leitner20a}. A comparison of a transformation- and a CQ-based BE method can be found in~\cite{schanz15a}.

The drawback of all BE formulations, either for elliptic and much stronger for hyperbolic problems, is the high storage and computing-time demand, as a standard formulation scales with $\mathcal{O}(M^2)$ for $M$ unknowns. In time domain, additionally, the time complexity has to be considered, where in the case of a CQ-based formulation, the complexity is of order $\mathcal{O}(M^2 N)$ for $N$ time steps. For elliptic problems, fast methods have been proposed, \eg the fast multipole method (FMM)~\cite{greengard97a}, or $\mathcal{H}$-matrix-based methods with the adaptive cross approximation (ACA) applied to the matrix blocks~\cite{bebendorf03,bebendorf08a}. The extension of FMM to the time variable has been published in~\cite{ergin98a} for acoustics and in~\cite{otani06a} for elastodynamics.
Fast methods for CQ are published in~\cite{banjai14a,messner10a}, involving a reformulation of CQ. 

In this paper, a different approach is used. Independently whether the CQ in its original form or gCQ is used, essentially, a three-dimensional data array needs to be computed and stored efficiently. This data array is determined by the spatial discretisation, resulting in two-dimensional data, and the selected complex frequencies, which gives the third dimension. To find a low-rank representation of this three-dimensional tensor, generalised adaptive cross approximation (3D-ACA) can be used. This technique is a generalisation 
of ACA~\cite{bebendorf03} and is proposed by Bebendorf et al.~\cite{bebndorf11a,bebendorf13a}.  It is based on a Tucker decomposition~\cite{tucker66a} and can be
traced back to the group of Tyrtyshnikov~\cite{oseledets08a}. The 3D-ACA can be seen as a higher-order singular value decomposition (SVD) or as a multilinear SVD~\cite{lathauwer00a}, which is a generalisation of the matrix SVD to tensors. Alternative but related approaches are presented in~\cite{dirckx22a} or~\cite{bauinger21a}, and do not employ a low-rank representation but interpolate with respect to the frequency to reduce the total number of frequency points needed for retrieving the time-domain solution. \corr{The approach of Anderson et al.~\cite{anderson20a} should be mentioned as well, where a frequency/time hybrid integral-equation method is presented.}

Here, the 3D-ACA is applied to a gCQ-based time-domain formulation utilising the original idea of the multivariate ACA~\cite{bebendorf13a}. The original version utilising ACA in the $\mathcal{H}$-matrices for the frequencies will be compared to a version where the FMM is used. Recently, a very similar approach has been published by Seibel~\cite{seibel22a}, where the conventional CQ method is used and, contrary to here, the $\mathcal{H}^2$ technique is used for approximating the BE matrices at the different frequencies. A short version using the indirect BE method can be found in~\cite{schawccm2024}. Here, a collocation approach based on the direct BE method is used for Dirichlet and mixed problems, whereas a Galerkin approach is used for Neumann problems. After briefly discussing the governing integral equations and their discretisation in space and time, the 3D-ACA is presented for both approaches. The main contribution of this paper is the comparison of both and the application to a real-world problem, \ie the sound scattering of an electrical machine. 
\section{Problem statement}
Let $\Omega \subset \mathbb{R}^3$ be a bounded Lipschitz domain and
$\Gamma = \partial\Omega$ its boundary with the outward normal
$\mathbf{n}$. The acoustic wave propagation is governed by
\begin{subequations}\label{eq:govtot}
  \begin{align} \label{eq:gov}
    \frac{\partial^2}{\partial t^2} u(\x,t)-
    c^2\Delta u(\x,t)& =0
                       \quad &(\x,t)&\in\Omega\times (0,T)\\
    u(\x,0)
    =\frac{\partial}{\partial t} u(\x,0)
                     &= 0
                       \quad & \x &\in \Omega \\
    \gamma_0 u(\x,t) &= g_D (\x,t) 
                       \quad & (\x,t) &\in \Gamma_D\times(0,T) \\
    \gamma_1 u(\x,t) &= g_N (\x,t) 
                       \quad & (\x,t) &\in \Gamma_N\times(0,T) 
  \end{align}
\end{subequations}
with the acoustic pressure $u(\x,t)$, the wave speed $c$ and the end time $T>0$. $\Gamma_D$ and
$\Gamma_N$ denote, respectively, the Dirichlet and Neumann
boundary, with $\Gamma = \Gamma_D \cup \Gamma_N$ and $\Gamma_D \cap
\Gamma_N =\emptyset$. The respective traces are defined by
\begin{align} \label{eq:traceD}
  \gamma_0 u(\x,t) = \lim_{\Omega \owns \x \to
  \x \in \Gamma} u(\x,t)
\end{align}
for the Dirichlet trace and by
\begin{align} \label{eq:traceN}
  \gamma_{1,x} u(\x,t) = \lim_{\Omega \owns \x \to
  \x \in \Gamma} \kl{\nabla u(\x,t) \cdot \mathbf{n}} = q(\x,t)
\end{align}
for the Neumann trace, also known as the conormal derivative or flux
$q(\x,t)$. Note, in the engineering literature the flux is usually defined with a negative sign in front of the nabla operator. Here, for consistency with the mathematical literature the definition without the negative sign is used. Hence, an outwards flux is positive.

Problem~\eqref{eq:govtot} can be solved with integral equations (see,
\eg~\cite{sauterschwab11, bonnet99a,sayas2016}), either by a direct approach or an indirect approach using layer potentials. Here, the former will be used. First, the retarded boundary integral operators are formulated. The single-layer potential is defined by
\begin{align} \label{eq:slp}
  (\op{V} * q)(\x,t) = 
  \int_0^t \int_\Gamma \gamma_0U(\x - \y, t - \tau) q(\y, \tau) \Od{s}_\y
  \Od{\tau} \; ,
\end{align}
with the fundamental solution $U(\x - \y, t - \tau) = \frac{1}{4 \pi \|\x-\y\|} \delta\kl{t-\tau -
  \frac{\|\x-\y\|}{c}}$, where $\delta(\cdot)$ denotes the Dirac distribution, which introduces a weak singularity. The double-layer potential is defined by
\begin{align} \label{eq:dlp}
  (\op{K} * u)(\x,t) = \int_0^t  \int_{\Gamma} \gamma_{1,y} U(\x - \y, t -
  \tau) u( \y, \tau) \Od{s}_\y \Od{\tau} 
\end{align}
and its adjoint
\begin{align} \label{eq:adlp}
  (\op{K}' * u)(\x,t) = \int_0^t  \gamma_{1,x}  \int_{\Gamma} U(\x - \y, t -
  \tau) u( \y, \tau) \Od{s}_\y \Od{\tau} \; .
\end{align}
The remaining operator to be defined is the hypersingular operator
\begin{align} \label{eq:hso}
  (\op{D} * u)(\x,t) = - \int_0^t  \gamma_{1,x}  \int_{\Gamma} \; \gamma_{1,y} U(\x - \y, t -
  \tau) u( \y, \tau) \Od{s}_\y \Od{\tau} \; .
\end{align}
Based on the representation formula
\begin{align} \label{eq:representationBIE}
   u\kl{\x,t} = \kl{\op{V} \ast q} \!\!\kl{\x,t}  -\kl{\op{K} \ast u}\!\!\kl{\x,t} 
  \quad        (\x,t)\in\Omega\times (0,T) \; ,
\end{align}
and by either a collocation or Galerkin formulation, the
respective boundary integral equations can be established. For a
collocation approach, the starting point is
\begin{equation} \label{eq:bie_collo}
  \op{C}u \kl{\x,t} =  \kl{\op{V} \ast q} \!\!\kl{\x,t}  -
  \kl{\op{K} \ast u}\!\!\kl{\x,t}
  \quad        (\x,t)\in\Gamma\times (0,T) 
\end{equation}
where the limiting process to the boundary results in the so-called
\corr{integral-free term $\op{C}$, see, \eg~\cite{mantic93} how the
integral-free term depends on the geometry and its computation.} The above collocation
approach will be used to study the Dirichlet and mixed problems.
In addition, a Neumann problem in an open domain will be studied, \ie
 a scattering problem. For that, a Galerkin
formulation is preferred, which can be formulated using the 
conormal derivative of \eqref{eq:representationBIE}. Applying suitable
traces results in
\begin{align} \label{eq:bie_galerkin}
    \kl{\op{D} \ast u}\!\!\kl{\x,t} = \frac{1}{2} q \kl{\x,t} - \kl{\op{K}' \ast q} \!\!\kl{\x,t} 
    \quad        (\x,t)\in\Gamma_N\times (0,T) \;.
\end{align}
Note, the integral-free term in a Galerkin formulation simplifies to $1/2$. The hypersingular operator in
\eqref{eq:bie_galerkin} is regularised to a weakly singular one by partial integration,
see, \eg~\cite{steinbach08a}. 
\section{Boundary element method: Discretisation}

\paragraph{Spatial discretisation}
The boundary $\Gamma$ is discretised resulting in an
approximation
\begin{equation}\label{eq:geomapprox}
  \Gamma^h = \bigcup_{e=1}^{E} \tau_e\,,
\end{equation}
which is the union of $E$ geometrical boundary elements $\tau_e$, here linear surface triangles. Finite element bases on boundaries $\Gamma_D$ and $\Gamma_N$ are
used to construct the approximation spaces
\begin{align} \label{eq:spaces}
    X_D
  = \operatorname{span}\{\varphi_1,\varphi_2,\dots,\varphi_{M_1} \}
  \quad \text{and} \quad
  X_N
  =  \operatorname{span}\{\psi_1,\psi_2,\dots,\psi_{M_2} \}.
\end{align}
The unknowns $u$ and  $q$ are approximated by 
linear combinations of functions in $X_D$ and $X_N$:
\begin{equation}\label{eq:ansatz_space}
  u^h = \sum_{\ell=1}^{M_1} u_\ell(t) \varphi_\ell(\y) 
 \quad \text{and} \quad
 q^h = \sum_{k=1}^{M_2} q_k(t) \psi_k(\y)  \; .
\end{equation}
Note, the coefficients $u_\ell(t)$ and $q_k(t)$ are still continuous functions of
time $t$. In the following, linear continuous shape functions
$\varphi_\ell$ and constant discontinuous shape functions $\psi_k$ will be chosen, in order to achieve a consistent overall discretisation, see~\cite{sayas2016} for further details.

Inserting $\varphi_\ell$ and $\psi_k$ in \eqref{eq:bie_collo}, and applying the
collocation method at the nodal collocation points $\x_i$
results in the semi-discrete equation system 
\begin{equation}  \label{eq:bie_semi_collo}
  \begin{bmatrix} \mat{V}_{DD}(t) & -\mat{K}_{DN}(t) \\ \mat{V}_{ND}(t) &
    -\kl{\mat{C}_{NN} + \mat{K}_{NN}(t)} \end{bmatrix} \ast
  \begin{bmatrix} \mat{q}_D^h(t) \\ \mat{u}_N^h (t) \end{bmatrix} =
  \begin{bmatrix} \mat{C}_{DD} + \mat{K}_{DD}(t) & -\mat{V}_{DN}(t) \\ \mat{K}_{ND}(t) &
    - \mat{V}_{NN}(t)\end{bmatrix} \ast
  \begin{bmatrix} \mat{g}_D^h(t) \\ \mat{g}_N^h (t) \end{bmatrix}
\end{equation}
with the semi-discrete matrix entries
\begin{align*}
  \mat{V}(t)[i,j] &= \int\limits_{\operatorname{supp} \psi_j}  \gamma_0 U(\x_i - \y, t) \psi_j(\y)   \Od{s}_\y \\
  \mat{K}(t)[i,j] &= \int\limits_{\operatorname{supp} \varphi_j}  \gamma_{1,y} U(\x_i - \y, t) \varphi_j(\y)   \Od{s}_\y \; ,
\end{align*}
where, here and in the following, the sans serif font denotes either
vectors or matrices collecting the nodal values.
Furthermore, the indices $D$ and $N$ distinguish between nodes at $\Gamma_D$ and $\Gamma_N$, respectively. Hence, in \eqref{eq:bie_semi_collo}, a splitting with respect
to known and unknown data corresponding to $\Gamma_D$ and
$\Gamma_N$ is carried out. Note, the convolution with
respect to time is still not discretised and the matrix entries
$\mat{V}(t)[i,j]$ and $\mat{K}(t)[i,j]$ are thus continuous functions of time.

In case of the Galerkin approach, $\varphi_\ell$ and $\psi_k$ are inserted in the integral equation \eqref{eq:bie_galerkin}. The resulting semi-discrete equation reads
\begin{equation}  \label{eq:bie_semi_galerkin}
  \mat{D}_{NN}(t)  \ast  \mat{u}_N^h (t) =
    \frac{1}{2}\mat{I} \mat{g}_N^h (t) - \mat{K}'_{NN}(t) \ast \mat{g}_N^h (t) \,,
\end{equation}
with the time-dependent matrix entries
\begin{align*}
    \mat{D}_{NN}(t)[i,j]   &= - \!\! \int\limits_{\operatorname{supp} \varphi_i}
    \varphi_i(\x) \gamma_{1,x} \!\! \int\limits_{\operatorname{supp} \varphi_j}
    \hspace{-4.65ex} = \;
    \gamma_{1,y} U(\x - \y, t) \varphi_j(\y)   \Od{s}_\y \Od{s}_\x\\ 
    \mat{K}_{NN}'(t)[i,j] &= \!\! \int\limits_{\operatorname{supp} \varphi_i}
                       \varphi_i(\x) \gamma_{1,x}\!\! \!\!\int\limits_{\operatorname{supp} \varphi_j}  \!U(\x - \y, t) \varphi_j(\y)   \Od{s}_\y  \Od{s}_\x \, .
\end{align*}
Note, the indices $NN$ in \eqref{eq:bie_semi_galerkin} are superfluous
as for a pure Neumann problem, the boundary is not subdivided. They are, however, kept for preserving consistency with \eqref{eq:bie_semi_collo}.

In the Galerkin formulation, double integrations are needed,
whereas in the collocation scheme, only one spatial integration is
necessary. In the collocation method, the weakly singular integrals are treated with a Duffy transformation~\cite{duffy82a}, whereas, in the Galerkin method, they are treated with the formula by Erichson and Sauter~\cite{erichsen98}. The regular integrals are calculated using a standard Gaussian quadrature using a heuristic distance-based formula to determine the number of Gauss points. No further treatment of quasi-singular integrals is considered.

\paragraph{Temporal discretisation}
The semi-discrete integral equations are discretised
in time using gCQ~\cite{lopezfernandez13a, lopezfernandez15b}. Here, the
variant using Runge-Kutta methods as the underlying time-stepping
technique is applied. Here, only a brief sketch of the algorithm is given. A more extension description using the same notation can be found in~\cite{schanz24a}. Methological details and a numerical analysis are provided in~\cite{lopezfernandez15b}.

The gCQ algorithm is exemplarily explained for the matrix vector product
\begin{equation} \label{eq:exampleConv}
   (\mat{V}_{DD} \ast \mat{q}^h_D)(t) = \int_0^t 
   \mat{V}_{DD}(t - \tau) \mat{q}^h_D(\tau) \Od{\tau} = \mat{f}(t) \, .
 \end{equation}
The matrix $\mat{V}_{DD}$ is of
size $M_2 \times M_2$. The vectors $\mat{q}^h_D$ and $\mat{f}$ are of size
$M_2$. Note, the right-hand-side vector $\mat{f}(t)$ is the result of
the integral or is a given right-hand side if $\mat{q}^h_D$ is the unknown.

Time is discretised in $N$, not necessarily constant time steps
$\Dt_i$, \ie 
\[ [0,T]=[0, t_1, t_2, \ldots, t_N], \quad  \Dt_i=t_i - t_{i-1}, \;\; 
  i=1,2, \ldots, N \, . \] 
Consider an A- and L-stable Runge-Kutta
method given by its Butcher tableau
\begin{tabular}{c|c}
  $\cvek$ & $\A$\\\hline    & $\bt$ 
\end{tabular}
with $\A \in \mathbb{R}^{m \times m}$, $\bvek,\cvek \in \mathbb{R}^m$
and $m$ the number of stages. The stability conditions require that $\bt \A^{-1} =
(0,0,\ldots,1)$ holds and imply $c_m=1$. These severe restrictions are partly prescribed by numerical
analysis, but are mostly based on 
experience~\cite{schanz01a}.

Some notation has to be introduced. A Laplace
transformed function is denoted with $\hat{()}$ and for the Laplace variable
holds $s \in \mathbb{C}, s.t. \Re{s} >0$. \Eg $\hat{\mat{V}}(s)$
means that the  fundamental solution $U(\x - \y, t - \tau)$ in the
integral above is replaced by $\hat{U}(\x - \y, s)$. 
Note, if the symbol $\hat{()}$ embraces the symbol of a matrix, it means that the Laplace
transform is applied on each of the entries, whereas if the argument is the matrix
$\kl{\Dt_n \A}^{-1}$, it means that in each entry, we get a matrix of
size $m \times m$. Further, the notation $(\mat{q}^h_D)_n$ denotes the
vector $\mat{q}^h_D$ at time $t_n$ and collects all nodal values at
all stages $m$. Considering this notation, the
whole algorithm can be given as in~\cite{lopezfernandez12a}:
\begin{itemize}
\item First step:
  \begin{equation*}
    (\mat{f})_1  = \hat{\mat{V}}_{DD}\kl{\kl{\Dt_1 \A}^{-1}} (\mat{q}^h_D)_1
  \end{equation*}
  with the implicit assumption of a zero initial condition.
\item For all steps $n=2,\ldots,N$, the algorithm has two substeps:
  \begin{enumerate}
  \item Update the solution vector $\mat{x}_{n-1}$ at all integration points $s_{\ell}$:
    \begin{equation} \label{eq:gCQ_xn}
      \mat{x}_{n-1}\kl{s_{\ell}} = \kl{\mat{I} - \Dt_{n-1} s_{\ell} \A}^{-1}
      \kl{\kl{\bt \A^{-1} \cdot \mat{x}_{n-2}(s_{\ell})} \bone
        + \Dt_{n-1} \A (\mat{q}^h_D)_{n-1} }
    \end{equation}
    for $\ell=1,\ldots,N_Q$ with the number of integration points
    $N_Q$ and the vector $\bone=(1,1,\ldots,1)^{\mathrm{T}}$ of size $m$. Note, with a slight abuse of notation, the vector
    $\mat{x}_{n-1}$ has to be understood as the collection of the
    results at all stages $m$ for each node in $\mat{q}^h_D$.
  \item Compute $(\mat{f})_n $ in case $\mat{q}^h_D$ is given, or solve
    the system in case $(\mat{f})_n $ is given:
    \begin{equation} \label{eq:gcq_formula}
        (\mat{f})_n =  \hat{\mat{V}}_{DD}\kl{\kl{\Dt_n \A}^{-1}}
        (\mat{q}^h_D)_n + \sum\limits_{\ell=1}^{N_Q}
        \hat{\mat{V}}_{DD}\kl{s_\ell}   \matWeight{\ell}{(\mat{q}^h_D)_{n-1}}
        \; .
      \end{equation}
      The abbreviation $\matWeight{\ell}{(\mat{q}^h_D)_{n-1}}=\omega_{\ell} \kl{\bt \A^{-1} \cdot
      \mat{x}_{n-1}(s_{\ell})} \kl{\mat{I} - \Dt_n s_{\ell} \A}^{-1}
    \bone$  denotes a vector of size $M_2 m$, which depends on
    $\mat{x}_{n-1}(s_{\ell})$. Hence, also here, the results at all
    stages $m$ are given for all nodes in $\mat{q}^h_D$.
  \end{enumerate}
\end{itemize}
Essentially, this algorithm requires the evaluation of the integral kernel at $N_Q$ points $s_\ell$,
representing complex frequencies. Consequently, we get an array of system matrices $\hat{\mat{V}}_{DD}\kl{s_\ell}$ of
size $M_2 \times M_2 \times N_Q$ in addition to the matrix $ \hat{\mat{V}}_{DD}\kl{\kl{\Dt_n \A}^{-1}}$.
This array of system matrices will be interpreted as a three-dimensional
array of data and approximated by a data-sparse representation based on 3D-ACA.

Before proceeding to the data-sparsification algorithm, the discrete sets of integral
equations are given explicitly. For the collocation approach
\eqref{eq:bie_semi_collo}, these are
\begin{equation}  \label{eq:bie_discr_collo}
  \begin{split}
    \begin{bmatrix} \hat{\mat{V}}_{DD} & -\hat{\mat{K}}_{DN} \\ \hat{\mat{V}}_{ND} &
      -\kl{\mat{C}_{NN} + \hat{\mat{K}}_{NN}} \end{bmatrix}\kl{\kl{\Dt_n \A}^{-1}}
    \begin{bmatrix} \mat{q}_D^h \\ \mat{u}_N^h \end{bmatrix}_n
    = \begin{bmatrix} \mat{C}_{DD} + \hat{\mat{K}}_{DD} & -\hat{\mat{V}}_{DN} \\ \hat{\mat{K}}_{ND} &
      - \hat{\mat{V}}_{NN}\end{bmatrix}\kl{\kl{\Dt_n \A}^{-1}}
    \begin{bmatrix} \mat{g}_D^h\\ \mat{g}_N^h \end{bmatrix}_n  \\+
    \sum\limits_{\ell=1}^{N_Q} \kl{
      \begin{bmatrix} \hat{\mat{K}}_{DD} & -\hat{\mat{V}}_{DN} \\ \hat{\mat{K}}_{ND} &
        -
        \hat{\mat{V}}_{NN}\end{bmatrix}\kl{s_\ell} \begin{bmatrix}\matWeight{\ell}{(\mat{g}_D)_{n-1}}
        \\ \matWeight{\ell}{(\mat{g}_N)_{n-1}} \end{bmatrix} -
      \begin{bmatrix} \hat{\mat{V}}_{DD} & -\hat{\mat{K}}_{DN} \\ \hat{\mat{V}}_{ND} &
        -\hat{\mat{K}}_{NN} \end{bmatrix}\kl{s_\ell} 
      \begin{bmatrix}\matWeight{\ell}{(\mat{q}^h)_{n-1}} \\
        \matWeight{\ell}{(\mat{u}^h)_{n-1}} \end{bmatrix}} \; ,
  \end{split}
\end{equation}
and for the Galerkin approach,
\begin{equation}  \label{eq:bie_discr_galerkin}
  \begin{split}
      \hat{\mat{D}}_{NN} \kl{\kl{\Dt_n \A}^{-1}} (\mat{u}_N^h)_n
    = &\kl{\frac{1}{2}\mat{I} - \hat{\mat{K}}'_{NN}}\kl{\kl{\Dt_n \A}^{-1}}
    (\mat{g}_N^h)_n  \\ &-
    \sum\limits_{\ell=1}^{N_Q} \kl{
        \hat{\mat{K}}'_{NN}\kl{s_\ell} \matWeight{\ell}{(\mat{g}_N)_{n-1}} +
       \hat{\mat{D}}_{NN} \kl{s_\ell}
        \matWeight{\ell}{(\mat{u}^h)_{n-1}} } \; .
  \end{split}
\end{equation}
By examination, the computational cost can be estimated. The computation
of the matrices is $\mathcal{O}((N_Q+1) M^2)$. The evaluation of the time-stepping
method is $\mathcal{O}(N_Q N)$ matrix-vector multiplications. Note, as given
in \autoref{app:gcq}, for Runge-Kutta methods with stages $m>1$, a suitable choice is $N_Q=N
(\log(N))^2$.  The solution of
\eqref{eq:bie_discr_collo} or \eqref{eq:bie_discr_galerkin} contributes with
$\mathcal{O}(M^2 n_\text{iter})$. An iterative solver is necessary as long as the
  matrix of the actual time step is approximated with the FMM or a
  $\mathcal{H}$-technique. In these estimates, the spatial
dimension is denoted by $M^2$, which is either $M_1^2, M_2^2$ or $M_1M_2$, according to the formulation. The value $n_\text{iter}$
is the number of iterations for the equation solver, which is
usually small. 
\section{Three-dimensional adaptive cross approximation}
A low-rank approximation of a three-dimensional array of data or a tensor of
third order $\mathcal{C}\in\mathbb{C}^{M\times M\times N_Q}$
has been proposed in~\cite{bebendorf13a}
and is referred to as a generalisation of adaptive cross approximation,
also called  3D-ACA. The 3D array of data is generated along an outer product, \ie
\begin{equation} \label{eq:outerProd}
  \mathcal{C} = \mathbf{H} \otimes \mathbf{f} \,,
\end{equation}
with $\mathbf{H}\in\mathbb{C}^{M\times M}$ and $\mathbf{f}\in
\mathbb{C}^{N_Q}$. The matrix $\mathbf{H}$ corresponds to the spatial
discretisation of the potentials in the BE formulations at a specific frequency $s_{\ell}$,
\eg the single-layer potential in \eqref{eq:gcq_formula} on the
right-hand side. This matrix will be called face or slice and may be computed as a
dense matrix or an approximation thereof. Here, two choices will be considered: 
\begin{itemize}
\item \textbf{The face is an $\mathcal{H}$-matrix}

  Basically, the concept is to decompose the matrix into subblocks
  first, and then perform a low-rank approximation to obtain suitable
  subblocks (see, \eg~\cite{bebendorf08a}). The matrix is partitioned by a recursive subdivision of the geometry. This
  implies a decomposition of the array of degrees of freedom based on
  a certain strategy. The index set $I_0$, \eg $I_0 = \{1, \ldots,
  N_p\}$, is subdivided into two son clusters based on principal
  component analysis (PCA). Recursively performing this procedure
  generates a balanced cluster tree. A son cluster is not further
  subdivided if it would drop below a prescribed minimal size $b_{\mathrm{min}}$. Then, it is
  referred to as a leaf. After creating the cluster tree, the block
  cluster tree or the hierarchical structure of the matrix is
  constructed with the aid of the distance criterion
  \begin{equation} \label{eq:admissible}
    \operatorname{min}\{\operatorname{diam}(Cl_x),
    \operatorname{diam}(Cl_y)\} \le \eta
    \operatorname{dist}(Cl_x,Cl_y) \; ,
  \end{equation}
  with a given parameter $\eta \in \mathbb{R}_+$. The corresponding
  index set selecting cluster $Cl_x$ is denoted by $I$ and selecting cluster $Cl_y$
  by $J$. In this way, the indices of the matrix are permuted such that a hierarchy of blocks arises, which are classified as
  far-field or near-field interactions. A subblock $\mathbf{H}_{I \times J}$ of $\mathbf{H}$
  contributes to the near-field interaction
  if \eqref{eq:admissible} is not fulfilled. Such matrix blocks are
  computed by the standard BE method and stored without
  approximation. A subblock $\mathbf{H}_{I \times J}$ represents far-field interactions if
  \eqref{eq:admissible} is fulfilled. Such matrix blocks are approximated
  by adaptive cross approximation (ACA)~\cite{bebendorf03} with
  recompression, leading to a representation in an SVD-like form $\mathbf{H}_{I \times J}  \approx
  \mathbf{H}_{I \times J}^r= \mathbf{U} \Sigma \mathbf{V}^H$. Starting from the
  low-rank representation of a block $\mathbf{H}_{I \times J}^r = \mathbf{U}
  \mathbf{V}$, a QR-decomposition of the low-rank matrix $\mathbf{U}$
  and $ \mathbf{V}$ is performed 
  \begin{align}
    \mathbf{H}_{I \times J}  \approx
    \mathbf{H}_{I \times J}^r  = \mathbf{U} \mathbf{V}^H=
    \mathbf{Q}_U( \mathbf{R}_U \mathbf{R}_V^H ) \mathbf{Q}_V^H 
    = \mathbf{Q}_U\,\check{\mathbf{U}}\check{\Sigma}\check{\mathbf{V}}^H \mathbf{Q}_V^H
    = \overline{\mathbf{U}}\check{\Sigma}\overline{\mathbf{V}}^H
  \end{align}
  and, secondly, the SVD is applied on the smaller inner matrix $ \mathbf{R}_U  \mathbf{R}_V^H$. In the remainder of this paper, this version will be called ''ACA-based 3D-ACA'' or in short ''3D-ACA''. In the examples, the parameters are chosen as $b_{\mathrm{min}}=20$ and $\eta=0.8$.
  
\item \textbf{The face is approximated by the fast multipole
    method (FMM)}

  In this version of the algorithm, the matrix $\mathbf{H}$ is not
  constructed. Instead, its matrix-vector product is approximated by a black-box FMM~\cite{fong09a}. As usual in FMM,  a uniform
  cluster tree based on the usual geometric subdivision is
  applied. Hence, as above, a block cluster tree is established and
  with this somehow an $\mathcal{H}$-matrix. The kernel expansion in
  the admissible blocks is performed with a
  Chebychev interpolation denoted with $S_p$ of polynomial order $p$, \ie
  \begin{align} \label{eq:cheb}
    \hat{U}\kl{\x,\y, s_\ell} \approx \sum_nS_p(\x, \x_n) \sum_m \hat{U}\kl{\x_n,\y_m, s_\ell}  S_p(\y, \y_m) \; ,
  \end{align}
  where $\x_n, \y_m$ are the Chebychev nodes. A multi-level schema is used for calculating the
  matrix-vector products. In this formulation,
  only the M2L-operator, which is essentially the middle term in
  \eqref{eq:cheb}, is frequency dependent because the clustering is
  only based on a geometrical argument similar to
  \eqref{eq:admissible}, and the interpolation is performed with the
  same order $p$ for all frequencies. Below, this version will be called ''FMM-based 3D-ACA''. When no confusion can arise, it will be called in short ''3D-FMM'', although strictly speaking the compression along the frequency dimension remains ACA.
\end{itemize}
The vector $\mathbf{f}$, called fiber in
the following, collects selected elements of $\mathbf{H}$ for the set
of frequencies determined by the gCQ. The 
latter would amount to $N_Q$ entries and $N_Q$
faces. On the one hand, the 3D-ACA approximates the faces by low-rank matrices (ACA-based 3D-ACA) or by the FMM (FMM-based 3D-ACA). More importantly, the number of necessary frequencies is determined adaptively. Hence, a sum $\mathcal{C}^{(\ell)}$ of an increasing number $\ell$ of outer-product terms as in \eqref{eq:outerProd} is constructed. The sum is truncated when $\mathcal{C}^{(\ell)}$ approximates
$\mathcal{C}$ up to a prescribed precision $\varepsilon$, as measured with a Frobenius norm. Assuming some monotonicity, the norm can be computed recursively~\cite{seibel22a}, \ie
\begin{equation} \label{eq:norm}
  \left\| \mathcal{C}^{(\ell)}\right\|_F^2 = \sum\limits_{d,d'}^\ell \kl{\sum\limits_{i,j} H_d[i,j]
    \widebar{H_{d'}[i,j]}} \kl{\sum\limits_k f_d[k] \widebar{f_{d'}[k]}} \; .
\end{equation}

The overall 3D-ACA approach is sketched
in \autoref{algo:3d_aca} \corr{and is discussed in detail in~\cite{seibel22a,schanz24a}}.
\begin{algorithm}
  \caption{\label{algo:3d_aca}Pseudo code of 3D-ACA (taken from~\cite{seibel22a})}
  \begin{algorithmic}
    \Function{3D ACA}{$ENTRY, \varepsilon$}       \Comment{$ENTRY$
      provides the integrated kernel values at collocation
      point $x_i$ and element $j$}
    \State $\mathcal{C}^{(0)} =0, k_1=0$ and $\ell=0$
    \While{$\|\mathbf{H}_{\ell}\|_F \|\mathbf{f}_{\ell}\|_2 >
      \varepsilon \|\mathcal{C}^{(\ell)}\|_F$}
    \State $\ell = \ell+1$
    \State $H_{\ell}[i,j] = ENTRY(i,j,k_{\ell}) - \mathcal{C}^{(\ell -
      1)}[i,j,k_{\ell}], \qquad i,j = 1, \ldots , M$
    \State $H_{\ell}[i_{\ell},j_{\ell}] = \max_{i,j} | H_{\ell}[i,j]|$
    \State $f_{\ell}[k] = H_{\ell}[i_{\ell},j_{\ell}]^{-1} \kl{ENTRY(i_{\ell},j_{\ell},k) - \mathcal{C}^{(\ell -
        1)}[i_{\ell},j_{\ell},k]}, \qquad k = 1, \ldots , N_Q$
    \State $\mathcal{C}^{(\ell)} = \mathcal{C}^{(\ell-1)} +
    \mathbf{H}_{\ell}  \otimes \mathbf{f}_{\ell}$
    \State $k_{\ell+1} = \arg \max_k |f_{\ell}[k]|$
    \EndWhile
    \State $r = \ell -1$   \Comment{Final rank, \ie necessary
      frequencies}
    
    \Return $\mathcal{C}^r = \sum\limits_{\ell=1}^r \mathbf{H}_{\ell}  \otimes \mathbf{f}_{\ell}$
    \EndFunction
  \end{algorithmic}
\end{algorithm}
The 3D-ACA is not applied on the whole matrix
but on the matrix blocks defined by the block cluster trees. Hence,
each block can have a different rank $r$, \ie a different number of
selected frequencies. For near-field interactions, the algorithm is straight
forward. However, for far-field interactions, the question arises how the pivot
element $H_{\ell}[i_{\ell}, j_{\ell}]$ can be found. In the FMM-based
variant, the only relevant entries are the M2L-operators and,
hence, a maximum entry can be found. However, in the ACA-based variant, there is
no clear rule how to determine the maximum entry in an 
approximated block of the $\mathcal{H}$-matrix. In~\cite{schanz24a}, an estimate based on
a column- or row-wise consideration is suggested. However,
numerical tests have shown that these estimates are not sharp
enough. Here, the respective matrix block is computed out of the
low-rank representation, the pivot element is determined, and
afterwards the dense block is discarded again.

The multiplication of the three-dimensional data array with a vector is carried out by
the proposed algorithm. Essentially, the algorithm separates the
frequency dependency such that $\mathbf{H}_\ell\kl{\hat{\mat{V}}}$ does not depend on the frequency, whereas the fibers
$\mathbf{f}_\ell\kl{\hat{\mat{V}}}$ explicitly introduce this dependency. Let us use the multiplication on the right-hand side
of \eqref{eq:gcq_formula} as example. The multiplication is changed into
\begin{align} \label{eq:product}
\sum\limits_{\ell=1}^{N_Q} \hat{\mat{V}}\kl{s_\ell}
  \matWeightO{\ell} =
  \sum\limits_{\ell=1}^{N_Q}  \sum\limits_{k=1}^r \mathbf{H}_k \kl{\hat{\mat{V}}}
  \otimes f_k[\ell]\!\!\kl{\hat{\mat{V}}}
  \matWeightO{\ell}
  = \sum\limits_{k=1}^r \mathbf{H}_k \kl{\hat{\mat{V}}}
  \otimes  \sum\limits_{\ell=1}^{N_Q} f_k[\ell]\!\!\kl{\hat{\mat{V}}}
  \matWeightO{\ell} \; .
      \end{align}
The complexity of the original
operation is $\mathcal{O}\kl{N_Q M^2}$ for $M$ spatial unknowns. The approximated
version has the complexity $\mathcal{O}\kl{r (M^2 + N_Q)}$. It
consists of the inner sum, which requires  $N_Q$ multiplications, and the outer sum, which is a matrix-times-vector multiplication of size $M$. Hence, the leading term with $M^2$ has
only a factor $r$ instead of $N_Q$ for dense
multiplication. For larger problem sizes, the significant reduction from $N_Q$ to $r$ complex frequencies offers a substantially faster discrete convolution.

\section{Numerical examples}
The acceleration of the gCQ-based time-domain
BEM is first tested on a unit-cube example, both for a pure Dirichlet problem and for a mixed boundary value
problem. The
main goal is to show that the approximation of the 3D-ACA does not
spoil the results, \ie the newly introduced approximation error is
smaller than the discretisation error of the dense BE formulation itself. The main focus is on
the overall obtained reduction in storage, \ie the compression rate defined as the ratio between the storage with compression and the dense storage. The second example is a scattering
problem with prescribed Neumann data, \ie the sound radiation of an electric machine. In
all tests, both proposed formulations, ACA-based 3D-ACA (in short: 3D-ACA) and FMM-based 3D-ACA (in short: 3D-FMM), are compared. Linear triangles are used for discretising the geometry. Linear continuous and constant discontinuous shape functions are applied for the Cauchy data. The 2-stage Radau IIA is used as time-stepping method. The final equation is solved with BiCGstab without preconditioner. Its precision is
set equal to the approximation of the ACA in the respective level. 

\subsection{Unit cube loaded by a smooth pulse}\label{sec:numeric_cube}
\begin{figure}[htb]
  \subcaptionbox{Unit cube (level 1, $h=\unit[0.5]{m}$)}{
    \includegraphics[scale=.23]{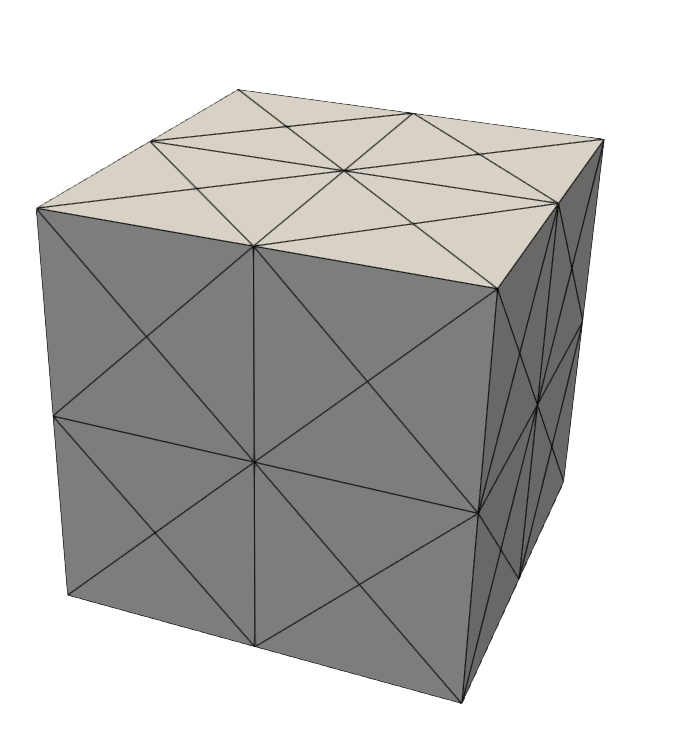}
  } 
  
  \subcaptionbox{\label{tab:refine} Used meshes}{
    \begin{tabular}{l|rrrrrr}
      \toprule
      refinement & nodes & elements & $h$ & $\Dt$ & $N$ & $N_Q$\\ \midrule
      1 & 50 & 96 & \unit[0.5]{m} & \unit[0.3]{s} & 10 & 27\\
      2 & 194 & 384 & \unit[0.25]{m} & \unit[0.15]{s} & 20 & 90 \\
      3 & 770 & 1536 & \unit[0.125]{m} & \unit[0.075]{s} & 40 & 272\\
      4 & 3074 & 6144 & \unit[0.0625]{m} & \unit[0.0375]{s} & 80 & 768\\
      5 & 12290 & 24576 & \unit[0.03125]{m} & \unit[0.01875]{s} & 160 & 2061\\
      \botrule
    \end{tabular}
  }
  \caption{\label{fig:unitCube}Unit cube: Geometry and discretisation parameters}
\end{figure}

The first test geometry is a unit cube $[-0.5, 0.5]^3$ centred around the origin. The wave speed
is set to $c=\unitfrac[1]{m}{s}$. The coarsest mesh is displaced in
\autoref{fig:unitCube}. The refinement levels and mesh data are given in the table aside. In each refinement 
level, the mesh size and the time step are halved. Hence, the ratio
of time step to mesh size was kept constant at $0.6$.
The meshes are created by bisecting the cathetus of the coarser
mesh. As load, a smooth pulse
\begin{align} \label{eq:load}
  u\kl{\mathbf{y},t} = \frac{\kl{t-\frac{r}{c}}^2}{r}  e^{-c
  \kl{t-\frac{r}{c}}} \; H\kl{t-\frac{r}{c}} \qquad \text{with} \quad
  r=\|\mathbf{x}-\mathbf{y}\| 
\end{align}
is exerted at the excitation point $\mathbf{x}=(0.8, 0.2, 0.3)^{\mathsf{T}}$ \corr{in order to dispose of an analytical solution for comparison}.  This load is applied on the whole boundary in case of the
Dirichlet problem and on one half of the
cube in case of the mixed problem. The other three sides are subjected to the conormal derivative
of \eqref{eq:load}. The total observation time is set to $T=\unit[3]{s}$ such
that the pulse travels over the whole unit cube.

To check that the compression does not deteriorate the accuracy of the results, the error in space
and time is measured with
\begin{align} \label{eq:error}
  \begin{split}
    L_{\mathrm{max}}(u_h)  &= \max\limits_{1 \le n \le N} \|
    u\kl{\frac{t_n+t_{n+1}}{2}} -  u_h\kl{\frac{t_n+t_{n+1}}{2}}
    \|_{L_2} \\ & = \max\limits_{1 \le n \le N} \sqrt{\int\limits_\Gamma 
      \kl{u\kl{\mathbf{x},\frac{t_n+t_{n-1}}{2}} - 
        u_h\kl{\mathbf{x},\frac{t_n+t_{n-1}}{2}} }^2 \Od{\Gamma} } \;
    ,
  \end{split}
\end{align}
where $u$ and $u_h$ are placeholders for the presented quantity
measured at the boundary of the selected boundary $\Gamma$.
The same error has been used in~\cite{sauterschanz17a}. It computes
the usual $L_2$-error with respect to the spatial variable and
selects the maximal value over all time steps.  The convergence
rate is denoted by 
\begin{equation} \label{eq:eoc}
  \mathrm{eoc} = \log_2\kl{\frac{L_{\mathrm{max}}^h}{L_{\mathrm{max}}^{h+1}} }
  \; ,
\end{equation}
where the superscripts $h$ and $h+1$ denote two subsequent refinement
levels. As long as the time response
is smooth, the convergence results 
presented below only change marginally 
if an $L_2$-error in time would be used. 

The approximation within the faces has been chosen to be
$\varepsilon_{\mathrm{ACA}}=10^{-4} \ldots 10^{-8}$ for the different
spatial refinement levels. The precision
of the method with respect to the frequencies, \ie the $\varepsilon$
in  Algorithm~\ref{algo:3d_aca} was selected as $\varepsilon =
100\varepsilon_{\mathrm{ACA}}$. In \autoref{tab:parameter}, the respective
parameters for each refinement level are given, together with the parameters of the FMM. Note, the indication \emph{levels} denotes the levels of the FMM multi-level schema and should not be confused with the refinement levels.
\begin{table}[hbt]
  \centering
  \begin{tabular}{l|cccc}
    \toprule
      refinement & $\varepsilon_{\mathrm{ACA}}$ & $\varepsilon$ & levels & order\\ \midrule
      1 & $10^{-4} $ & $10^{-2} $ & 1 & 1 \\
      2 & $10^{-5} $ & $10^{-3} $ & 2 & 2 \\
      3 & $10^{-6} $ & $10^{-4} $ & 2 & 3 \\
      4 & $10^{-7} $ & $10^{-5} $ & 3 & 4 \\
    5 & $10^{-8} $ & $10^{-6} $ & 4 & 5 \\
    \botrule
    \end{tabular}
  \caption{\label{tab:parameter}Selected parameters of the 3D-ACA/3D-FMM in
    each refinement level}
\end{table}

\begin{figure}[htb]
  \centering
  \subcaptionbox{Dirichlet problem}{
    \includegraphics{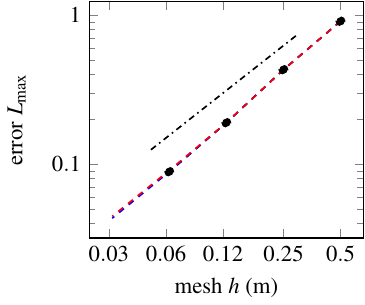}
  }
  \subcaptionbox{Mixed problem}{
    \includegraphics{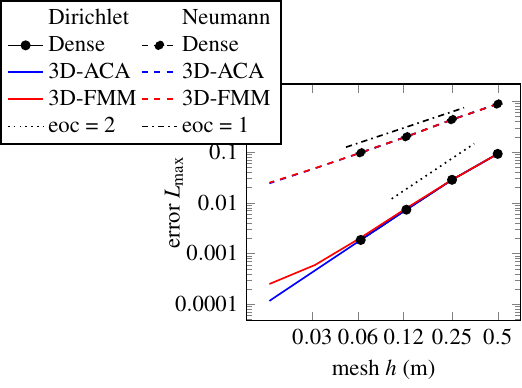}
  }
  \caption{Cube: $L_{\mathrm{max}}$-error versus refinement in space and time}
  \label{fig:ConvCube}
\end{figure}

The error defined in \eqref{eq:error} is plotted over the refined
mesh in \autoref{fig:ConvCube}. The spatial mesh size $h$ is displayed on the abscissa. The time-step size is reduced proportially (see~\autoref{tab:refine}).
The error \eqref{eq:error} measures a combined error is space and time. Hence, although the underlying Runge-Kutta method would allow a convergence order of
three, the dominating convergence orders are linear for
the Neumann data and quadratic for the Dirichlet data, according to the chosen low-order spatial shape functions. The numerical test shows that the 3D-ACA and 3D-FMM
do not spoil the overall accuracy. The lines indeed overlap with the
respective dense calculation.

\begin{figure}[htb]
  \subcaptionbox{Dirichlet problem}{
    \includegraphics{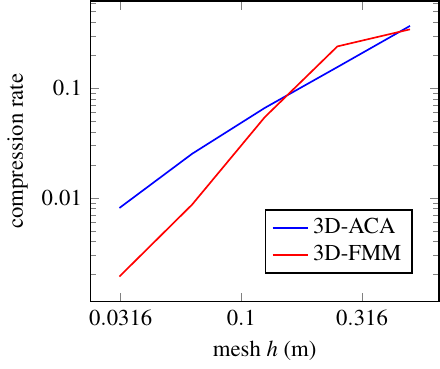}
  } 
  \subcaptionbox{Mixed problem}{
    \includegraphics{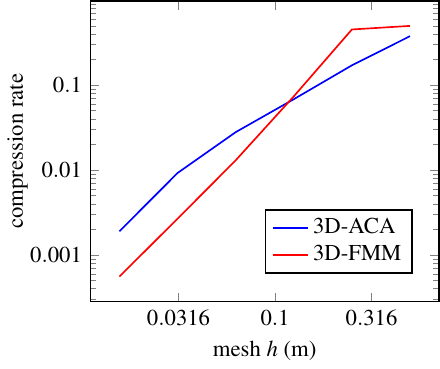}
    }
  \caption{\label{fig:CompressCube}Compression rates of 3D-ACA and 3D-FMM
    for both problems}
\end{figure}

\begin{figure}[htb]
  \centering
  \subcaptionbox{Single-layer potential} {
    \includegraphics{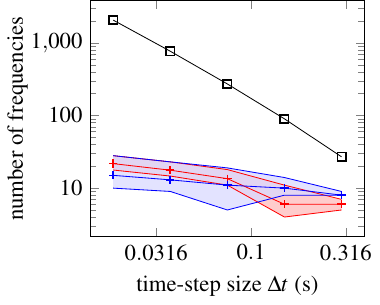}
  }
  \subcaptionbox{\label{fig:FreqsDiriDLP} Double-layer potential}{
    \includegraphics{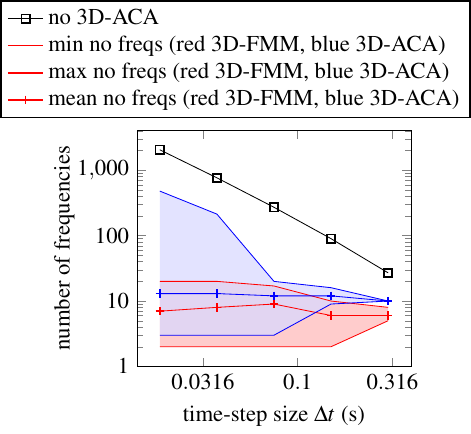}
  }
  \caption{Cube: Necessary number of frequencies for both approaches
    for the Dirichlet problem}
  \label{fig:FreqsDiri}
\end{figure}

The next question deals with the efficiency
of the 3D-ACA and 3D-FMM. In \autoref{fig:CompressCube}, the
compression rate is plotted versus the refinement.
As mentioned above, the compression rate is defined as the ratio of the memory used by 3D-ACA or 3D-FMM to the memory used by a dense computation. It is obvious that very high
compression rates can be obtained independently whether a Dirichlet
problem or a mixed problem is considered. The 3D-FMM results in better
compression rates if a sufficiently large problem is under study. \corr{This
  corresponds to our observations comparing ACA and FMM for elliptic problems.} The major contribution to the
compression rates results from the rank reduction with respect to the
frequencies. In \autoref{fig:FreqsDiri}, the number of
used/necessary frequencies is presented.
Note, as the algorithm is applied on each matrix block individually, the blocks have
different ranks. That is why in \autoref{fig:FreqsDiri}, a maximal, minimal,
and mean value is plotted. It is obvious that the mean value is very
small for both operators. However, the maximal value of the 3D-ACA for the
double-layer potential sometimes gets very large. This is caused
by the error computation \eqref{eq:norm}, which is critical if the
matrix entries are small or partially zero. This effect is not visible
for the 3D-FMM because only the M2L-operators are
considered for far-field interactions and they are the same as for the single-layer operator. This effect results from the implementation, where the
FMM of the double-layer potential is computed by shifting the normal
derivative to the interpolation polynomials
(see~\cite{schanz18a}). Hence, in \autoref{fig:FreqsDiriDLP}, only the
numbers for the near-field interactions of the 3D-FMM are displayed. Overall, the
mean value of the rank shows only a very slight increase with
increasing problem size. This very positive observation has also be
reported for different geometries in~\cite{schanz24a}. For the mixed problems, the plots show
exactly the same tendency, and are omitted here. It can be
concluded that the 3D-ACA or 3D-FMM are very effective in determining
the necessary number of frequencies.

\begin{figure}[htb]
  \centering
  \subcaptionbox{3D-FMM SLP}{
    \includegraphics{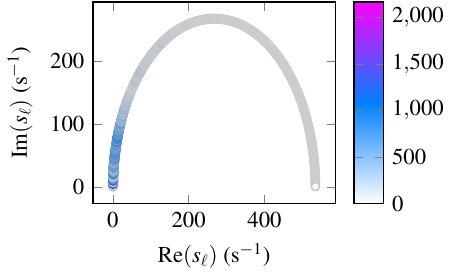}
  }
  \subcaptionbox{3D-ACA SLP}{
    \includegraphics{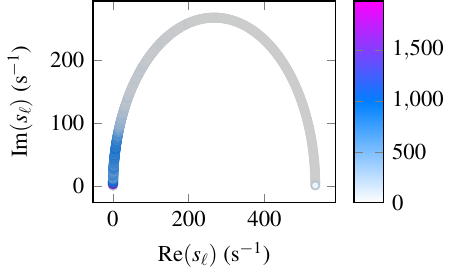}
  }
  \subcaptionbox{3D-FMM DLP}{
    \includegraphics{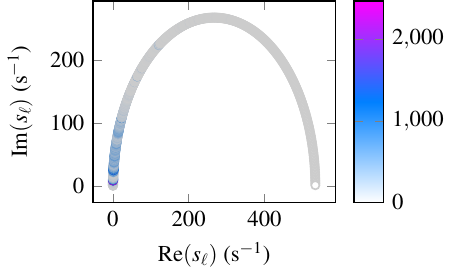}
  }
  \subcaptionbox{3D-ACA DLP}{
    \includegraphics{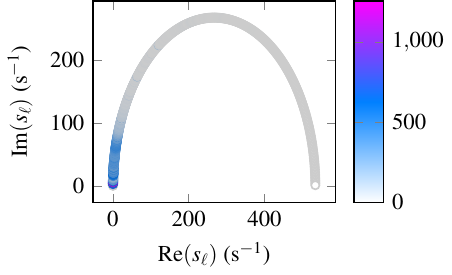}
  }
  \caption{\label{fig:CubeFreqsUsed}Cube (level 4): Used complex frequencies.
    The colour code corresponds to the number of matrix blocks at which the frequency is active: Dirichlet problem, single-layer potential (SLP, top) and double-layer potential (DLP, bottom), 3D-FMM (left) and 3D-ACA (right)}
\end{figure}

The reduction of the number of required frequencies raises the question which
frequencies are selected. The algorithm selects
the frequencies in an adaptive way like a greedy algorithm. The frequencies in the gCQ lie on a circle in the complex half plane, of which radius and location are determined by the parameters given in
\autoref{app:gcq}. Only half of the circles have to be determined
because the complex frequencies appear as pairs of complex
conjugates. In \autoref{fig:CubeFreqsUsed}, the half circles are shown with a colour coding indicating the number of matrix blocks that select the respective frequency.
Note, the first and second complex frequency are excluded as they are
selected from almost all blocks and would spoil the colour scale. It
is obvious that mostly frequencies with a small real  
part are selected. 3D-ACA or 3D-FMM select
similar frequencies. Moreover, there are only marginal differences
between the single-layer and the double-layer potential formulation.

\begin{figure}[hbt]
  \subcaptionbox{Dirichlet problem}{
    \includegraphics{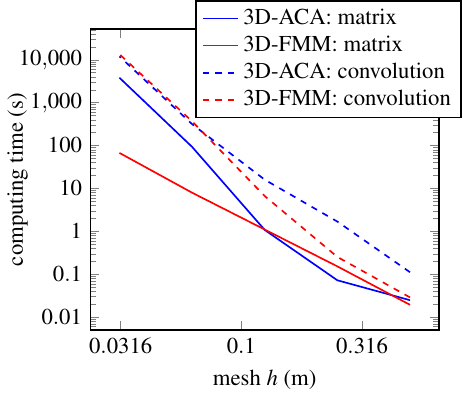} 
  } 
  \subcaptionbox{Mixed problem}{
    \includegraphics{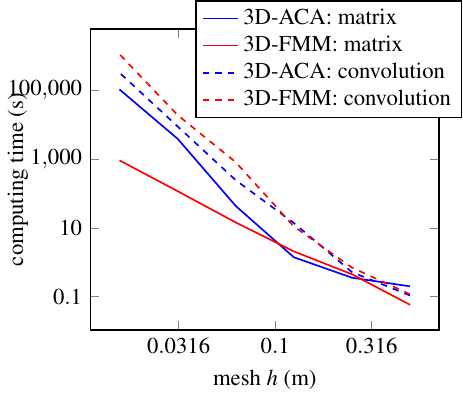} 
    }
  \caption{\label{fig:CPUTimeCube}Computing time for assembling and
    solving}
\end{figure}

The last aspect of the comparison concerns the
required computing time. In \autoref{fig:CPUTimeCube}, the computing time is plotted
against the spatial mesh size.
Essentially, the two major contributions within the algorithm are
compared. First, the time to establish the data-sparse representation
of the 3D-data array is studied, denoted with `matrix'. Obviously, the FMM-based version is faster than the
ACA-based version. This is in line with the behaviour of ACA and FMM for
elliptic problems. Essentially, the 3D-ACA has to establish and fill
the $\mathcal{H}$-matrices and the fibers, whereas the 3D-FMM has only
to establish the set of M2L operators and the fibers. The second step
is the convolution, where the tensor product has to be
computed as discussed in \eqref{eq:product}. This step is faster
for the 3D-ACA than for the 3D-FMM. This behaviour has as
well its analogue in the elliptic case. Nevertheless, this step
consumes the most computing time for larger problems and seems to be a
bottleneck of the 3D-FMM approach. The actual implementation of the FMM
is very standard and uses a uniform cluster tree, where an expert
selection of the levels is necessary and most probably
sub-optimal. Hence, this behaviour may be improved by a different
parameter set or, preferably, by a balanced cluster tree.

\subsection{Sound scattering of an electric machine}

\begin{figure}
  \centering
    \subcaptionbox{Geometry}{
      \includegraphics[trim={4cm 2cm 8cm 1cm},clip, scale=.3]{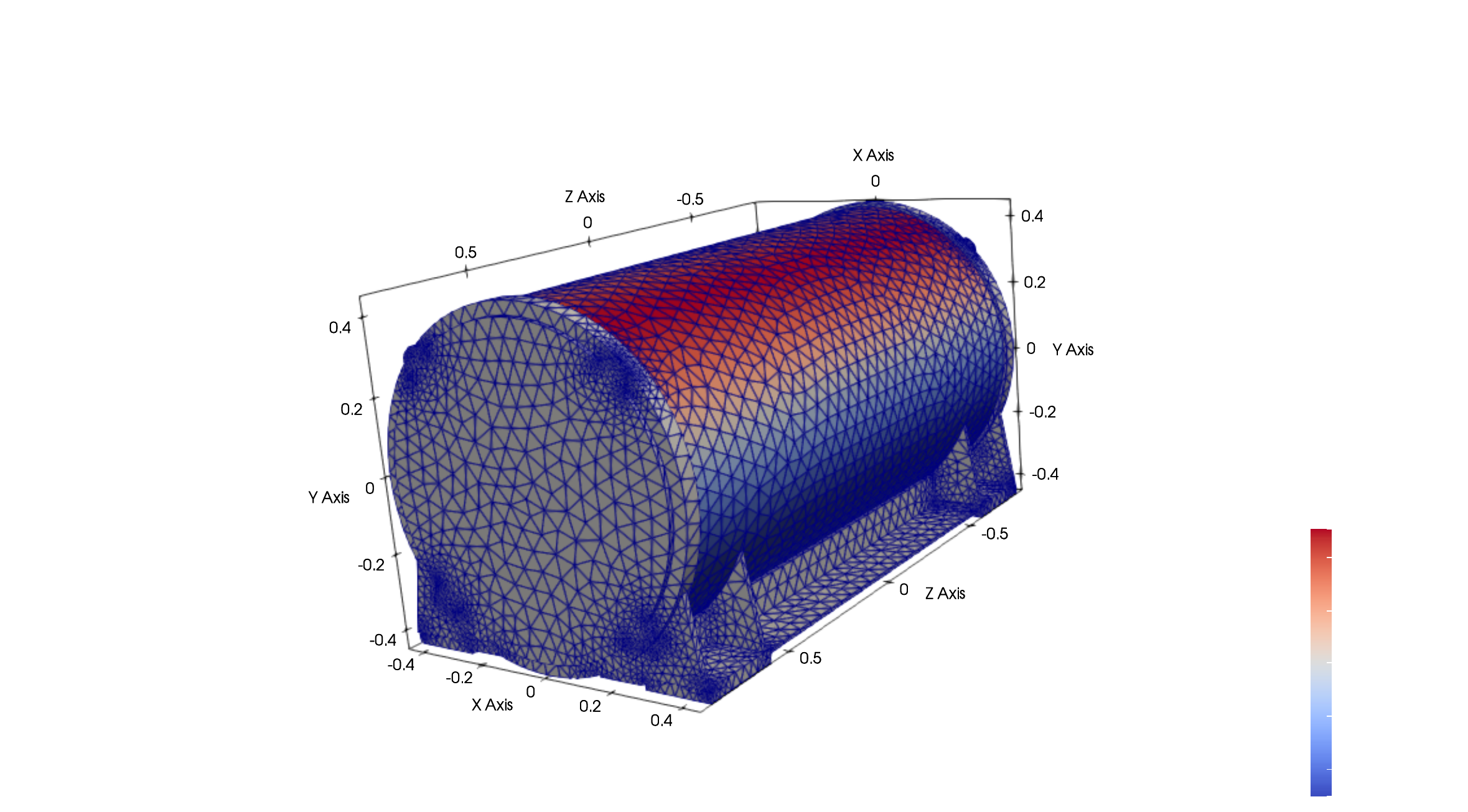}
  } 
  \subcaptionbox{\label{tab:mach_data}Computing data}{
    \begin{tabular}{ll}
      nodes & 49464 \\
      elements & 16488 \\
      mean $h$ & $\approx$ \unit[0.06]{m} \\
      $\Dt$ & \unit[0.0001]{s} \\
      $N$ & 150 \\
      $\varepsilon_{\mathrm{ACA}}$ & $10^{-5}$ \\
      $\varepsilon$ & $10^{-3}$ \\
      level, order & 4, 4
    \end{tabular}
  }
  \caption{\label{fig:machine}Electric machine: Dimensions,
    coordinate system and parameters}
\end{figure}

Next, a scattering example is
presented. The scatterer is an electric machine,
similar to the one in~\cite{weilharter11a} where the structural
vibrations and the 
acoustic scattering have been studied in frequency domain 
with a coupled FEM-BEM approach. \corr{Because the feed-back of the acoustic field to the structural vibration of the machine is negligible, both simulations steps can be carried out subsequently. Here, only the time-domain calculation of the acoustic scattering is performed. The excitation by the vibration of the surface carried over by a prescribed flux boundary condition}. The model is presented in
\autoref{fig:machine} together with its dimensions and the
data for the mesh and time stepping.
The wave speed of air
$c=\unitfrac[343]{m}{s}$ is assumed. The observation time is set to
$T=\unit[0.015]{s}$.

\begin{figure}
  \centering
  \subcaptionbox{\label{fig:fluxBC}Surface flux}{
    \includegraphics[trim={8cm 1cm 10cm 1cm},clip,
    scale=.35]{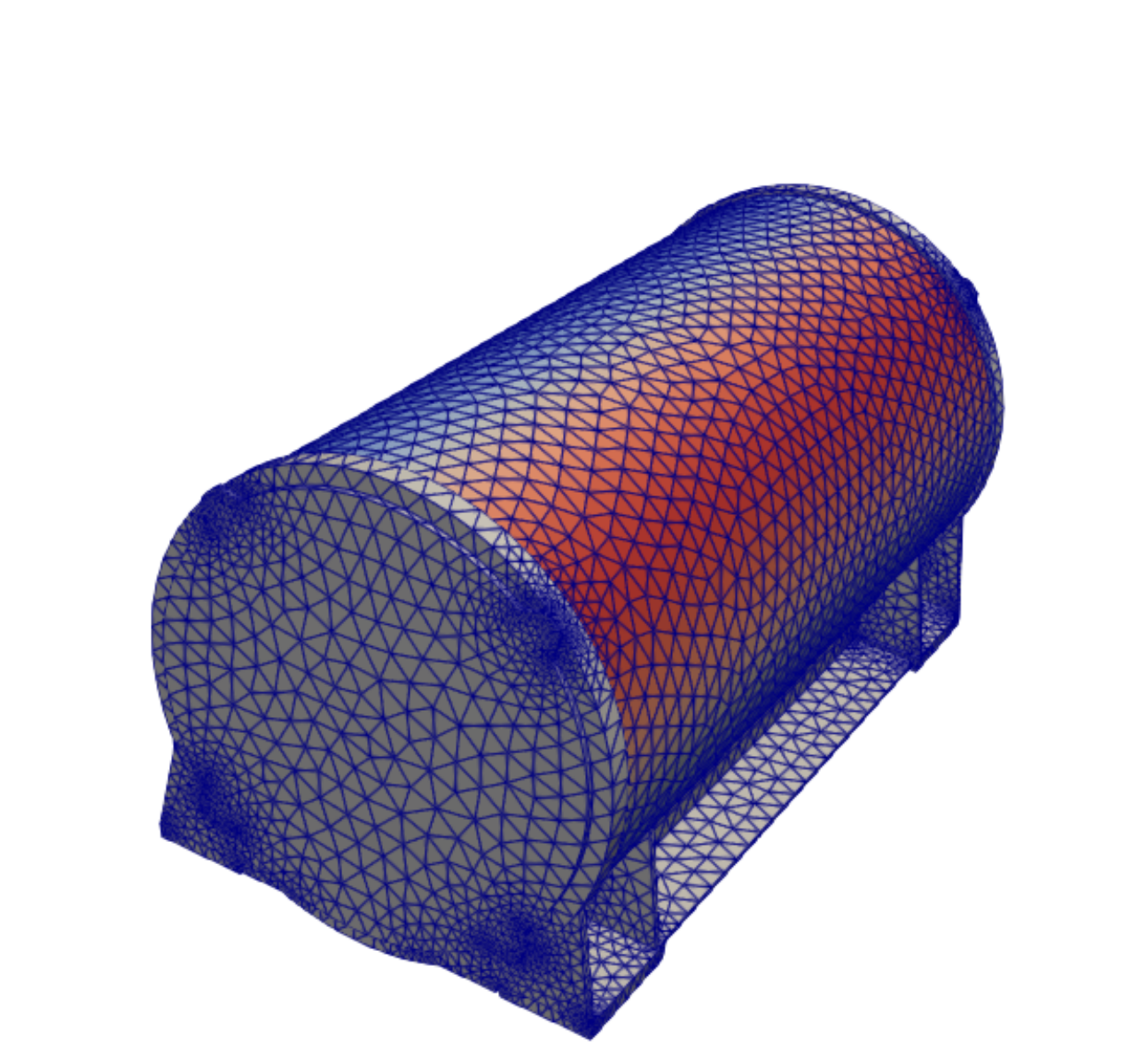}
  } 
  \subcaptionbox{\label{fig:timehistoryBC}Dependency on time}{
    \includegraphics{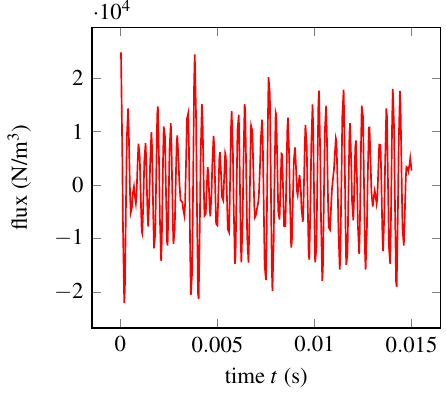}
  }
  \caption{\label{fig:machine_excite}Electric machine: Boundary
    conditions in space and time}
\end{figure}

The excitation is a radial quadrupole mode, which often appears as a
dominating vibration mode in an electric machine. The surface
displacement in normal direction is prescribed by 
\begin{equation} \label{eq:quad_modes}
  u_\text{surface} = \sum\limits_{i=1}^4 \sin(2 \pi \omega_i t - 2 \varphi) 40
  \cdot \unit[10^{-6}]{m} \quad
  \text{with} \;\; \omega_i \in \{2637.04, 2093.04, 2349.28, 3136.00\}
  [s^{-1}] 
\end{equation}
on the surface of the cylinder, whereas the other parts of the
geometry, \ie the end shields, cooling ribs and mounting feet, are
assumed to be rigid.  The circumferential angle $\varphi$ refers to a cylindrical
coordinate system aligned with the machine. Using the Euler equation $\frac{\partial
  u}{\partial n} = \varrho_\text{f} \ddot{u}_{\mathrm{surface}}$ with $\varrho_\text{f}$ the mass density of air, the displacement of the machine's surface can be carried over to a Neumann boundary condition for flux, \ie
\begin{equation} \label{eq:machine_bc}
  q(\varphi,t) = - \varrho_\text{f} \sum\limits_{i=1}^4 \kl{2 \pi\omega_i}^2 \sin(2 \pi \omega_i t - 2
  \varphi)  40 \cdot \unit[10^{-6}]{m}
\end{equation}
on the cylindrical hull and zero elsewhere. Note, the notation above with the pressure
denoted by $u$ and the surface displacement denoted by
$u_\text{surface}$ may cause a confusion. The excitation by the flux
is displayed in \autoref{fig:machine_excite}.
Note, the norms of the flux are shown for the midpoint of the observation period. The colour coding shows that
only the cylindrical surface is loaded by \eqref{eq:machine_bc}.

\begin{figure}
  \centering
   \subcaptionbox{Time $t=\unit[0.0011]{s}$}{
    \includegraphics[trim={6cm 1cm 0cm 1cm},clip,
    scale=.5]{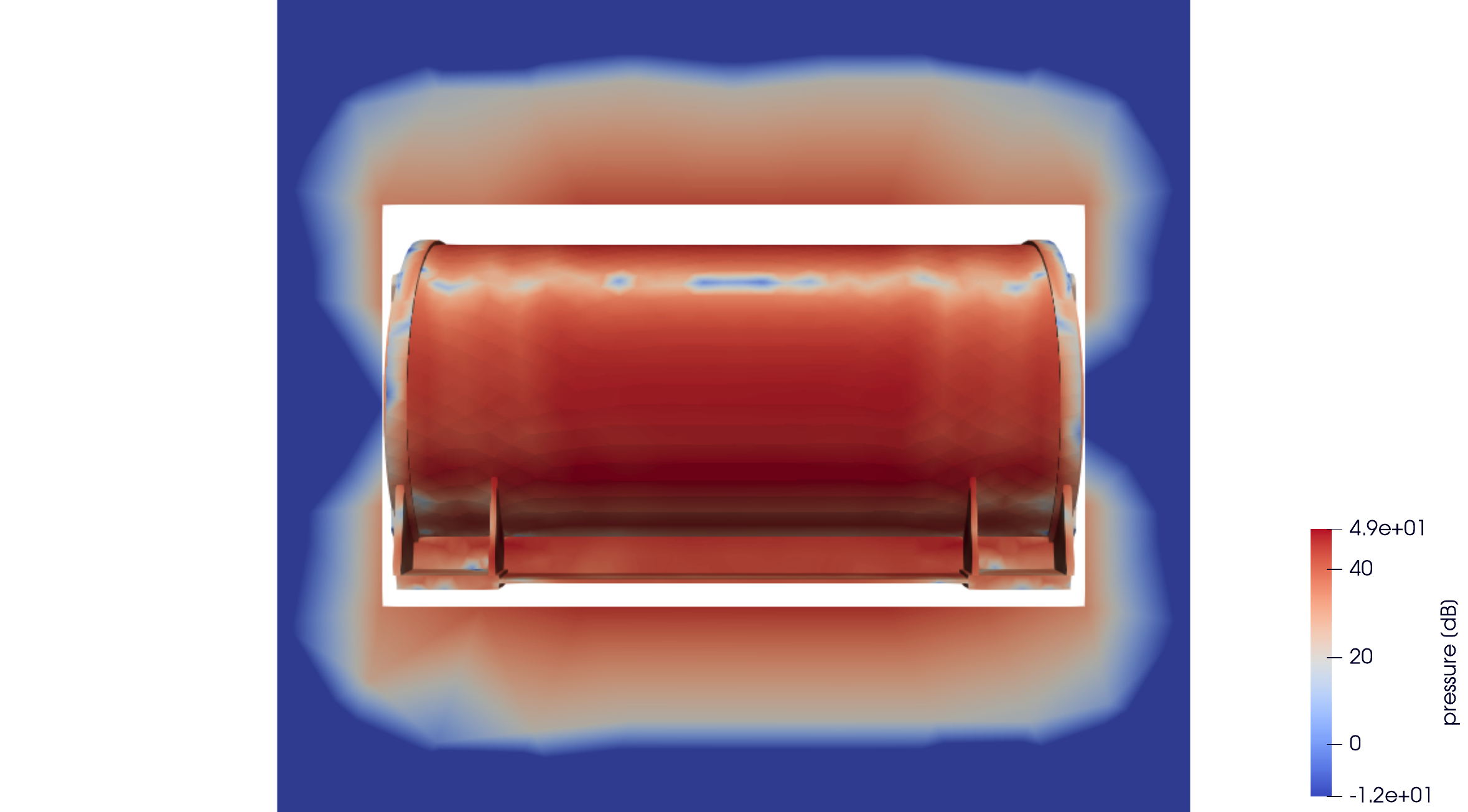}
    }
  \subcaptionbox{Time $t=\unit[0.0124]{s}$}{
    \includegraphics[trim={6cm 1cm 0cm 1cm},clip,
    scale=.5]{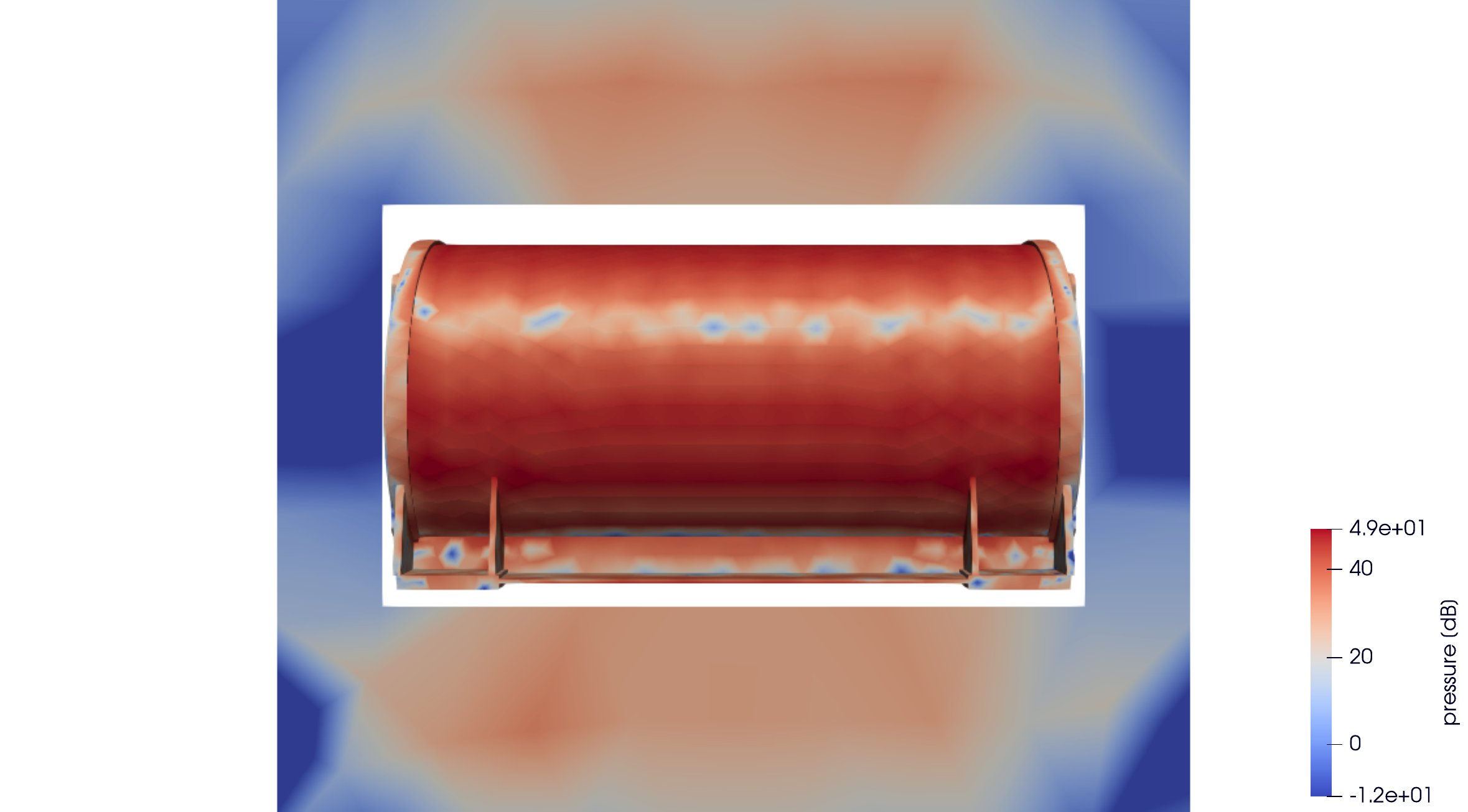}
    }
  \caption{Pressure (in $\unit{dB}$) on the surface of the electric machine and in the
    surrounding air (3D-ACA)}
  \label{fig:SoundPressSpace}
\end{figure}

First, results for a computation with 3D-ACA are presented in
\autoref{fig:SoundPressSpace} (For a movie see~\cite{schanzvideo25}). The sound pressure level is displayed on the
machine surface and in the surrounding air. The usual definition of the sound
pressure level $L_p = \unit[20
\log_{10}\kl{\frac{u}{p_0}}]{dB}$ with $p_0=\unitfrac[2 \cdot
10^{-5}]{N}{m^2}$ is used.
The snapshots show the pressure field at different time instants.
As expected, the values at the surface are larger than the scattered ones. Note, the scattered values have been determined by evaluating the
representation formula \eqref{eq:representationBIE} after computing
the Cauchy data, \ie the pressure field on the surface. The partly  visible
lines are attributed to the visualisation programme, which interpolates
the pressure values between the selected points linearly. Nevertheless, the scattered field looks reasonable. At the first selected time instant, the pressure wave has not yet travelled through the whole domain, \ie a wave front can be observed. The second selected time instant shows the scattered field after the wave has travelled over the whole computed domain and thus shows an established sound field. The results obtained by
3D-FMM are equal and are omitted here.

\begin{figure}
  \centering
  \subcaptionbox{\label{fig:SoundPressTimeSurf}Point on the cylinder (node 3014)}{%
    \includegraphics{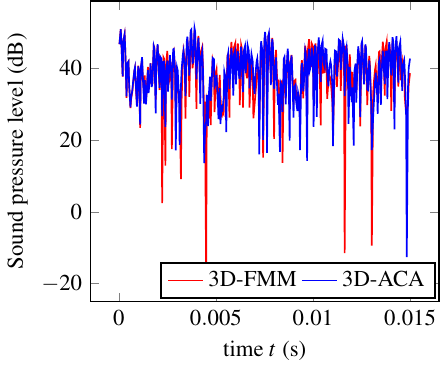}
  }
  \subcaptionbox{\label{fig:SoundPressTimeInnen}Point in the domain
    (node 66)}{
    \includegraphics{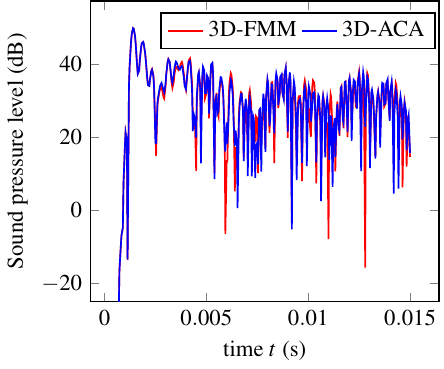}
  }
  \caption{Sound pressure level (in $\unit{dB}$) over time}
  \label{fig:SoundPressTime}
\end{figure}

\begin{figure}
  \centering
  \includegraphics{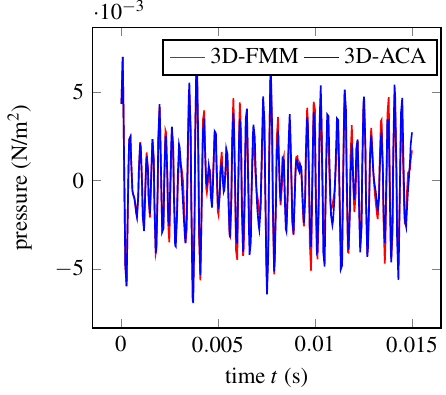}
  \caption{Pressure on the surface over time}
  \label{fig:PressTime}
\end{figure}

Next, the sound pressure level over time is presented in
\autoref{fig:SoundPressTime}.
Two immediate observations can be made. First, 3D-ACA
and 3D-FMM seem to provide different sound pressures. This is
true, although it must be remarked that the logarithmic scale amplifies the differences.
\corr{To better show the differences, the
  surface pressure at the same point (node 3014) as in
  \autoref{fig:SoundPressTimeSurf} is shown in \autoref{fig:PressTime}.}
 \corr{It can be observed that the pressure values deviate in the peak
  values only marginally. This is explained by the fact that 3D-ACA and 3D-FMM introduce different approximations. More precisely, the necessary rank
  with respect to the complex frequencies is determined based on the
  norm \eqref{eq:norm}, for which both methods measure different
  quantities (see the discussion above \eqref{eq:norm}). Several combinations of the 
  numerical parameters have been tested and led to similar results. In all tests, the overall accuracy as measured by \eqref{eq:norm} was met for both approaches.}
Furthermore, the results in
\autoref{fig:SoundPressTimeInnen} start from infinite values. This is
a correct result for a time-domain calculation. The signal needs some
time to propagate. This
causal behaviour results in a zero pressure for the first few time steps,
which corresponds to an infinite sound pressure level. As well, the sinusoidal
shape of the excitation is visible, where again at the
points in air, first, a kind of impulse is visible, which is related
to the arriving wave.

The 3D-ACA obtained a compression rate of 0.00878453 and needed
\unit[50503]{s} of computation. The 3D-FMM computation reached a compression
rate of 0.0105453 and needed \unit[176406]{s} of computation. This seems
to be contradictory to the results of
\autoref{sec:numeric_cube}. However, there, the parameters of
both approximations were set such that the same error is obtained, whereas here, some heuristic parameters are chosen, which give
comparable but not identical results. Especially, in the 3D-ACA, a very
weak precision has been selected, whereas in the 3D-FMM, a reliable
interpolation order has been selected. The latter is very conservative and
does not result in the best compression. However, as we do not have an
exact solution, the 3D-FMM with the selected
parameters is expected to be closer to the true solution than the ACA with the weak precision.

This electric-machine example illustrates that time-domain BEM is feasible in a large-scale technically relevant setting, in particular when compression by an ACA-based or FMM-based 3D-ACA is employed.
\section{Conclusions}
The complexity of the time-domain BEM is
$\mathcal{O}(M^2 N)$ for $N$ time steps and $M$ spatial degrees of
freedom. The generalised adaptive cross approximation (3D-ACA) has been applied to obtain a data-sparse and fast method. In addition to a compression in space by a standard ACA or fast multipole method (FMM), the 3D-ACA reduces the set of applied complex frequencies of the gCQ adaptively, leading to a nearly linear complexity in space and
time. The comparison of the FMM- and ACA-based versions show a clear advantage of the FMM-based one with respect to storage, or of the ACA-based one with respect to computation time. A further improvement of the FMM-based version may include an adaptive tree or a different form of the FMM. The calculation of the sound scattered by an
electric machine proves that 3D-ACA enables the application of time-domain BEM to large-scale examples.

\paragraph*{Acknowledgement}
This work is supported by the joint DFG/FWF Collaborative Research Centre CREATOR (DFG: Project-ID 492661287/TRR 361; FWF: 10.55776/F90) at TU~Darmstadt, TU~Graz and JKU~Linz.

\appendix
\section{Parameters of the gCQ} \label{app:gcq}
The derivation and reasoning how the integration weights and points
are determined can be found in~\cite{lopezfernandez12a, lopezfernandez15b}. The result of
these papers are recalled here.
The integration points in the complex plane are
\begin{align*}
  s_\ell = \gamma(\sigma_\ell)\ , & & \omega_\ell = \frac{4 K(k^2)}{1
                                      \pi i} \gamma^\prime
                                      (\sigma_\ell)\ , && N_Q = N \log(N)\ , 
\end{align*}
where for Runge-Kutta methods with  $m>1$  stages, it should be $N_Q = N
(\log(N))^2$. $K(k)$ is the complete elliptic integral of first kind, \ie
\begin{align*}
 K(k) = \int^1_0 \frac{dx}{\sqrt{(1-x^2) (1-k^2 x^2)}}\ , & & K^\prime (k) = K(1-k)\ ,
\end{align*}
and $K^\prime$ is its derivative, which equals the integral of the complementary modulus.
The argument $k$ depends on the relation $q$ of the maximum and minimum step sizes in the following way:
\begin{align*}
 k = \frac{q - \sqrt{2q-1}}{q + \sqrt{2q-1}}  \qquad q = 5\, 
  \frac{\Dt_{\text{max}} \max_i
  |\lambda_i\kl{\A}|}{\Dt_{\text{min}}\min_i
  |\lambda_i\kl{\A}|} \; ,
\end{align*}
with the eigenvalues $\lambda_i\kl{\A}$. For the implicit Euler method,
the eigenvalues are 1 and the factor 5 in $q$ can be skipped. 
The functions $\gamma(\sigma_\ell)$ and $\gamma^\prime(\sigma_\ell)$ are
\begin{align*}
 \gamma(\sigma_\ell) =& \frac{1}{\Dt_{\text{min}}(q-1)} \left(
                        \sqrt{2q-1} \frac{k^{-1} +
                        \text{sn}(\sigma_\ell)}{k^{-1}-\text{sn}(\sigma_\ell)}
                        -1 \right)  \\[1ex]
 \gamma^\prime(\sigma_\ell) =&
                               \frac{\sqrt{2q-1}}{\Dt_{\text{min}}(q-1)}
                               \frac{2\ \text{cn}(\sigma_\ell)\
                               \text{dn}(\sigma_\ell)}{k (k^{-1}-
                               \text{sn}(\sigma_\ell))^2}  \\[1ex]
 \sigma_\ell =& -K(k^2) + \left(\ell-\frac{1}{2} \right) \frac{4
                K(k^2)}{N_Q} + \frac{i}{2} K^\prime(k^2) \, ,
\end{align*}
where $\text{sn}(\sigma_\ell)$, $\text{cn}(\sigma_\ell)$ and
$\text{dn}(\sigma_\ell)$ are the Jacobi elliptic functions. 
As seen above, the integration contour is only determined by the
largest and smallest time steps chosen, but does not depend on any
intermediate step sizes. Due to the symmetric distribution of the integration points with
respect to the real axis, only half of the frequencies $s_\ell$ need to be
calculated.

Last, it may be remarked that for constant time steps $\Dt_i=const$
and $\max_i |\lambda_i\kl{\A} = \min_i |\lambda_i\kl{\A}$, the parameter determination would fail
because this choice results in $q=1$ and $k=0$. This happens, \eg for
BDF1 with a constant time step size. Unfortunately, this
value is not allowed for the complete elliptic integral. However, a
slight change in the parameter $q \approx 1$ fixes this problem
without spoiling the algorithm. The latter can be done as these
parameter choices are empirical.

\bibliography{literature}

\begin{thebibliography}{47}
\providecommand{\natexlab}[1]{#1}
\providecommand{\url}[1]{\texttt{#1}}
\expandafter\ifx\csname urlstyle\endcsname\relax
  \providecommand{\doi}[1]{doi: #1}\else
  \providecommand{\doi}{doi: \begingroup \urlstyle{rm}\Url}\fi

\bibitem[Aimi and Diligenti(2008)]{aimi08a}
A.~Aimi and M.~Diligenti.
\newblock A new space-time energetic formulation for wave propagation analysis
  in layered media by {BEMs}.
\newblock \emph{Int. J. Numer. Methods. Engrg.}, 75\penalty0 (9):\penalty0
  1102--1132, 2008.

\bibitem[Aimi et~al.(2012)Aimi, Diligenti, Frangi, and Guardasoni]{aimi12a}
A.~Aimi, M.~Diligenti, A.~Frangi, and C.~Guardasoni.
\newblock A stable 3d energetic {Galerkin BEM} approach for wave propagation
  interior problems.
\newblock \emph{Eng. Anal. Bound. Elem.}, 36\penalty0 (12):\penalty0
  1756--1765, 2012.
\newblock ISSN 0955-7997.
\newblock \doi{http://dx.doi.org/10.1016/j.enganabound.2012.06.003}.
\newblock URL \url{http://www.sciencedirect.com/science/article/pii/
  S0955799712001361}.

\bibitem[Anderson et~al.(2020)Anderson, Bruno, and Lyon]{anderson20a}
Thomas~G. Anderson, Oscar~P. Bruno, and Mark Lyon.
\newblock High-order, dispersionless ``fast-hybrid'' wave equation solver. part
  i: O(1) sampling cost via incident-field windowing and recentering.
\newblock \emph{SIAM J. Sci. Comput.}, 42\penalty0 (2):\penalty0 A1348--A1379,
  2020.
\newblock \doi{10.1137/19M1251953}.
\newblock URL \url{https://doi.org/10.1137/19M1251953}.

\bibitem[Bamberger and Ha-{D}uong(1986)]{bambduong}
A.~Bamberger and T.~Ha-{D}uong.
\newblock Formulation variationelle espace-temps pour le calcul par potentiel
  retard{\'{e}} d'une onde acoustique.
\newblock \emph{Math. Meth. Appl. Sci.}, 8:\penalty0 405--435 and 598--608,
  1986.

\bibitem[Banjai and Kachanovska(2014)]{banjai14a}
L.~Banjai and M.~Kachanovska.
\newblock Fast convolution quadrature for the wave equation in three
  dimensions.
\newblock \emph{J. Comput. Phys.}, 279:\penalty0 103--126, 2014.
\newblock ISSN 0021-9991.
\newblock \doi{https://doi.org/10.1016/j.jcp.2014.08.049}.
\newblock URL \url{https://www.sciencedirect.com/science/article/pii/
  S0021999114006251}.

\bibitem[Bauinger and Bruno(2021)]{bauinger21a}
Christoph Bauinger and Oscar~P. Bruno.
\newblock ``interpolated factored green function'' method for accelerated
  solution of scattering problems.
\newblock \emph{J. Comput. Phys.}, 430:\penalty0 110095, 2021.
\newblock ISSN 0021-9991.
\newblock \doi{https://doi.org/10.1016/j.jcp.2020.110095}.
\newblock URL \url{https://www.sciencedirect.com/science/article/pii/
  S002199912030869X}.

\bibitem[Bebendorf(2008)]{bebendorf08a}
M.~Bebendorf.
\newblock \emph{Hierarchical Matrices: {A} Means to Efficiently Solve Elliptic
  Boundary Value Problems}, volume~63 of \emph{Lecture Notes in Computational
  Science and Engineering}.
\newblock Springer-Verlag, 2008.

\bibitem[Bebendorf(2011)]{bebndorf11a}
M.~Bebendorf.
\newblock Adaptive cross approximation of multivariate functions.
\newblock \emph{Constr. Approx.}, 34\penalty0 (2):\penalty0 149--179, 2011.
\newblock \doi{10.1007/s00365-010-9103-x}.
\newblock URL \url{https://doi.org/10.1007/s00365-010-9103-x}.

\bibitem[Bebendorf and Rjasanow(2003)]{bebendorf03}
M.~Bebendorf and S.~Rjasanow.
\newblock Adaptive low-rank approximation of collocation matrices.
\newblock \emph{Computing}, 70:\penalty0 1--24, 2003.

\bibitem[Bebendorf et~al.(2013)Bebendorf, K{\"u}hnemund, and
  Rjasanow]{bebendorf13a}
M.~Bebendorf, A.~K{\"u}hnemund, and S.~Rjasanow.
\newblock An equi-directional generalization of adaptive cross approximation
  for higher-order tensors.
\newblock \emph{Appl. Num. Math.}, 74:\penalty0 1--16, 2013.
\newblock ISSN 0168-9274.
\newblock \doi{https://doi.org/10.1016/j.apnum.2013.08.001}.
\newblock URL \url{http://www.sciencedirect.com/science/article/pii/
  S0168927413000950}.

\bibitem[Bonnet(1999)]{bonnet99a}
M.~Bonnet.
\newblock \emph{Boundary Integral Equation Methods for Solids and Fluids}.
\newblock John Wiley \& Sons, 1999.

\bibitem[Costabel(2005)]{costabel04}
M.~Costabel.
\newblock Time-dependent problems with the boundary integral equation method.
\newblock In E.~Stein, R.~de~Borst, and T.~J.~R. Hughes, editors,
  \emph{Encyclopedia of Computational Mechanics}, volume 1, Fundamentals,
  chapter~25, pages 703--721. John Wiley \& Sons, New York, Chichester,
  Weinheim, 2005.

\bibitem[Cruse and Rizzo(1968)]{cruri}
T.~A. Cruse and F.~J. Rizzo.
\newblock A direct formulation and numerical solution of the general transient
  elastodynamic problem, {I}.
\newblock \emph{Aust. J. Math. Anal. Appl.}, 22\penalty0 (1):\penalty0
  244--259, 1968.

\bibitem[De~Lathauwer et~al.(2000)De~Lathauwer, De~Moor, and
  Vandewalle]{lathauwer00a}
Lieven De~Lathauwer, Bart De~Moor, and Joos Vandewalle.
\newblock A multilinear singular value decomposition.
\newblock \emph{SIAM J. Matrix Aanal. A.}, 21\penalty0 (4):\penalty0
  1253--1278, 2000.
\newblock \doi{10.1137/S0895479896305696}.
\newblock URL \url{https://doi.org/10.1137/S0895479896305696}.

\bibitem[Dirckx et~al.(2022)Dirckx, Huybrechs, and Meerbergen]{dirckx22a}
Simon Dirckx, Daan Huybrechs, and Karl Meerbergen.
\newblock Frequency extraction for bem matrices arising from the 3d scalar
  helmholtz equation.
\newblock \emph{SIAM J. Sci. Comput.}, 44\penalty0 (5):\penalty0 B1282--B1311,
  2022.
\newblock \doi{10.1137/20M1382957}.
\newblock URL \url{https://doi.org/10.1137/20M1382957}.

\bibitem[Duffy(1982)]{duffy82a}
M.~G. Duffy.
\newblock Quadrature over a pyramid or cube of integrands with a singularity at
  a vertex.
\newblock \emph{SIAM J. Numer. Anal.}, 19\penalty0 (6):\penalty0 1260--1262,
  1982.

\bibitem[Ergin et~al.(1998)Ergin, Shanker, and Michielssen]{ergin98a}
A.~A. Ergin, B.~Shanker, and E.~Michielssen.
\newblock Fast evaluation of three-dimensional transient wave fields using
  diagonal translation operators.
\newblock \emph{J. Comput. Phys.}, 146\penalty0 (1):\penalty0 157--180, 1998.
\newblock \doi{10.1006/jcph.1998.5908}.

\bibitem[Erichsen and Sauter(1998)]{erichsen98}
S.~Erichsen and S.~A. Sauter.
\newblock Efficient automatic quadrature in 3-d {G}alerkin {BEM}.
\newblock \emph{Comput. Methods Appl. Mech. Engrg.}, 157\penalty0
  (3--4):\penalty0 215--224, 1998.

\bibitem[Fong and Darve(2009)]{fong09a}
W.~Fong and E.~Darve.
\newblock The black-box fast multipole method.
\newblock \emph{J. Comput. Phys.}, 228\penalty0 (23):\penalty0 8712--8725,
  2009.

\bibitem[Greengard and Rokhlin(1997)]{greengard97a}
L.~Greengard and V.~Rokhlin.
\newblock A new version of the {Fast Multipole Method} for the {L}aplace
  equation in three dimensions.
\newblock \emph{Acta Num.}, 6:\penalty0 229--269, 1997.

\bibitem[Haider et~al.(2024)Haider, Rjasanow, and Schanz]{schanz24a}
Anita~M. Haider, Sergej Rjasanow, and Martin Schanz.
\newblock Generalised adaptive cross approximation for convolution quadrature
  based boundary element formulation.
\newblock \emph{Comput. Math. Appl.}, 175:\penalty0 470--486, 2024.
\newblock ISSN 0898-1221.
\newblock \doi{https://doi.org/10.1016/j.camwa.2024.10.025}.
\newblock URL \url{https://www.sciencedirect.com/science/article/pii/
  S0898122124004681}.

\bibitem[Kupradze et~al.(1979)Kupradze, Gegelia, Basheleishvili, and
  Burchuladze]{kupradze79}
V.~D. Kupradze, T.~G. Gegelia, M.~O. Basheleishvili, and T.~V. Burchuladze.
\newblock \emph{Three-Dimensional Problems of the Mathematical Theory of
  Elasticity and Thermoelasticity}, volume~25 of \emph{Applied Mathematics and
  Mechanics}.
\newblock North-Holland, Amsterdam New York Oxford, 1979.

\bibitem[Leitner and Schanz(2021)]{leitner20a}
Michael Leitner and Martin Schanz.
\newblock Generalized convolution quadrature based boundary element method for
  uncoupled thermoelasticity.
\newblock \emph{Mech Syst Signal Pr}, 150:\penalty0 107234, 2021.
\newblock ISSN 0888-3270.
\newblock \doi{https://doi.org/10.1016/j.ymssp.2020.107234}.
\newblock URL \url{http://www.sciencedirect.com/science/article/pii/
  S0888327020306208}.

\bibitem[Lopez-Fernandez and Sauter(2013)]{lopezfernandez13a}
M.~Lopez-Fernandez and S.~Sauter.
\newblock Generalized convolution quadrature with variable time stepping.
\newblock \emph{IMA J. of Numer. Anal.}, 33\penalty0 (4):\penalty0 1156--1175,
  2013.
\newblock \doi{10.1093/imanum/drs034}.

\bibitem[Lopez-Fernandez and Sauter(2015)]{lopezfernandez12a}
M.~Lopez-Fernandez and S.~Sauter.
\newblock Generalized convolution quadrature with variable time stepping. part
  {II}: Algorithm and numerical results.
\newblock \emph{Appl. Num. Math.}, 94:\penalty0 88--105, 2015.

\bibitem[L{\'o}pez-Fern{\'a}ndez and Sauter(2016)]{lopezfernandez15b}
Mar{\'\i}a L{\'o}pez-Fern{\'a}ndez and Stefan Sauter.
\newblock Generalized convolution quadrature based on {Runge-Kutta} methods.
\newblock \emph{Numer. Math.}, 133\penalty0 (4):\penalty0 743--779, 2016.
\newblock \doi{10.1007/s00211-015-0761-2}.

\bibitem[Lubich(1988{\natexlab{a}})]{lubich88a}
C.~Lubich.
\newblock Convolution quadrature and discretized operational calculus. {I.}
\newblock \emph{Numer. Math.}, 52\penalty0 (2):\penalty0 129--145,
  1988{\natexlab{a}}.

\bibitem[Lubich(1988{\natexlab{b}})]{lubich88b}
C.~Lubich.
\newblock Convolution quadrature and discretized operational calculus. {II.}
\newblock \emph{Numer. Math.}, 52\penalty0 (4):\penalty0 413--425,
  1988{\natexlab{b}}.

\bibitem[Mansur(1983)]{man}
W.~J. Mansur.
\newblock \emph{A Time-Stepping Technique to Solve Wave Propagation Problems
  Using the Boundary Element Method}.
\newblock Phd thesis, University of Southampton, 1983.

\bibitem[Manti\v{c}(1993)]{mantic93}
V.~Manti\v{c}.
\newblock A new formula for the {C-matrix} in the somigliana identity.
\newblock \emph{J. Elasticity}, 33:\penalty0 191--201, 1993.

\bibitem[Messner and Schanz(2010)]{messner10a}
M.~Messner and M.~Schanz.
\newblock An accelerated symmetric time-domain boundary element formulation for
  elasticity.
\newblock \emph{Eng. Anal. Bound. Elem.}, 34\penalty0 (11):\penalty0 944--955,
  2010.
\newblock \doi{10.1016/j.enganabound.2010.06.007}.

\bibitem[Oseledets et~al.(2008)Oseledets, Savostianov, and
  Tyrtyshnikov]{oseledets08a}
I.~V. Oseledets, D.~V. Savostianov, and E.~E. Tyrtyshnikov.
\newblock Tucker dimensionality reduction of three-dimensional arrays in linear
  time.
\newblock \emph{SIAM J. Matrix Aanal. A.}, 30\penalty0 (3):\penalty0 939--956,
  2008.
\newblock \doi{10.1137/060655894}.
\newblock URL \url{https://doi.org/10.1137/060655894}.

\bibitem[Otani et~al.(2007)Otani, Takahashi, and Nishimura]{otani06a}
Y.~Otani, T.~Takahashi, and N.~Nishimura.
\newblock A fast boundary integral equation method for elastodynamics in time
  domain and its parallelisation.
\newblock In M.~Schanz and O.~Steinbach, editors, \emph{Boundary Element
  Analysis: Mathematical Aspects and Applications}, volume~29 of \emph{Lecture
  Notes in Applied and Computational Mechanics}, pages 161--185.
  Springer-Verlag, Berlin Heidelberg, 2007.

\bibitem[Peirce and Siebrits(1997)]{peirce97}
A.~Peirce and E.~Siebrits.
\newblock Stability analysis and design of time-stepping schemes for general
  elastodynamic boundary element models.
\newblock \emph{Int. J. Numer. Methods. Engrg.}, 40\penalty0 (2):\penalty0
  319--342, 1997.
\newblock
  \doi{10.1002/(SICI)1097-0207(19970130)40:2\%3C319::AID-NME67\%3E3.0.CO;2-I}.

\bibitem[Sauter and Schanz(2017)]{sauterschanz17a}
S.A. Sauter and M.~Schanz.
\newblock Convolution quadrature for the wave equation with impedance boundary
  conditions.
\newblock \emph{J. Comput. Phys.}, 334:\penalty0 442--459, 2017.
\newblock ISSN 0021-9991.
\newblock \doi{http://dx.doi.org/10.1016/j.jcp.2017.01.013}.
\newblock URL \url{//www.sciencedirect.com/science/article/pii/
  S0021999117300232}.

\bibitem[Sauter and Schwab(2011)]{sauterschwab11}
Stefan Sauter and Christoph Schwab.
\newblock \emph{Boundary Element Methods}.
\newblock Number~39 in Springer Series in Computational Mathematics. Springer
  Verlag, Heidelberg, 2011.
\newblock \doi{10.1007/978-3-540-68093-2}.

\bibitem[Sayas(2016)]{sayas2016}
F.-J. Sayas.
\newblock \emph{Retarded Potentials and Time Domain Boundary Integral
  Equations: A Road Map}, volume~50 of \emph{Springer Series in Computational
  Mathematics}.
\newblock Springer, Cham, 2016.
\newblock \doi{10.1007/978-3-319-26645-9}.

\bibitem[Schanz(2001)]{schanz01a}
M.~Schanz.
\newblock \emph{Wave Propagation in Viscoelastic and Poroelastic Continua: A
  Boundary Element Approach}, volume~2 of \emph{Lecture Notes in Applied
  Mechanics}.
\newblock Springer-Verlag, Berlin, Heidelberg, New York, 2001.
\newblock \doi{10.1007/978-3-540-44575-3}.

\bibitem[Schanz(2018)]{schanz18a}
M.~Schanz.
\newblock Fast multipole method for poroelastodynamics.
\newblock \emph{Eng. Anal. Bound. Elem.}, 89:\penalty0 50--59, 2018.
\newblock \doi{10.1016/j.enganabound.2018.01.014}.

\bibitem[Schanz and Antes(1997)]{schanz97e}
M.~Schanz and H.~Antes.
\newblock A new visco- and elastodynamic time domain boundary element
  formulation.
\newblock \emph{Comput. Mech.}, 20\penalty0 (5):\penalty0 452--459, 1997.
\newblock \doi{10.1007/s004660050265}.

\bibitem[Schanz and Keshava(2025)]{schanzvideo25}
M.~Schanz and V.~Lakshmi Keshava.
\newblock Sound radiation of an electric machine.
\newblock online available, 2025.
\newblock URL \url{https://doi.org/10.3217/hhp1q-efw32}.

\bibitem[Schanz(2024)]{schawccm2024}
Martin Schanz.
\newblock {3D-ACA} for the time domain boundary element method: Comparison of
  {FMM} and $\mathcal{H}$-matrix based approaches.
\newblock In \emph{WCCM2024}, 2024.
\newblock URL \url{https://www.scipedia.com/public/_2024a}.

\bibitem[Schanz et~al.(2015)Schanz, Ye, and Xiao]{schanz15a}
Martin Schanz, Wenjing Ye, and Jinyou Xiao.
\newblock Comparison of the convolution quadrature method and enhanced inverse
  {FFT} with application in elastodynamic boundary element method.
\newblock \emph{Comput. Mech.}, 57\penalty0 (4):\penalty0 523--536, 2015.
\newblock ISSN 0178-7675.
\newblock \doi{10.1007/s00466-015-1237-z}.
\newblock URL \url{http://dx.doi.org/10.1007/s00466-015-1237-z}.

\bibitem[Seibel(2022)]{seibel22a}
Daniel Seibel.
\newblock Boundary element methods for the wave equation based on hierarchical
  matrices and adaptive cross approximation.
\newblock \emph{Numer. Math.}, 150\penalty0 (2):\penalty0 629--670, 2022.
\newblock \doi{10.1007/s00211-021-01259-8}.
\newblock URL \url{https://doi.org/10.1007/s00211-021-01259-8}.

\bibitem[Steinbach(2008)]{steinbach08a}
O.~Steinbach.
\newblock \emph{Numerical Approximation Methods for Elliptic Boundary Value
  Problems}.
\newblock Springer, 2008.

\bibitem[Tucker(1966)]{tucker66a}
Ledyard~R. Tucker.
\newblock Some mathematical notes on three-mode factor analysis.
\newblock \emph{Psychometrika}, 31\penalty0 (3):\penalty0 279--311, 1966.
\newblock \doi{10.1007/BF02289464}.
\newblock URL \url{https://doi.org/10.1007/BF02289464}.

\bibitem[Weilharter et~al.(2011)Weilharter, B{\'\i}r{\'o}, Lang, and
  Rainer]{weilharter11a}
Bernhard Weilharter, Oszk{\'a}r B{\'\i}r{\'o}, Hermann Lang, and Siegfried
  Rainer.
\newblock Computation of the noise radiation of an induction machine using {3D
  FEM/BEM}.
\newblock \emph{COMPEL - The international journal for computation and
  mathematics in electrical and electronic engineering}, 30\penalty0
  (6):\penalty0 1737--1750, 2025/02/20 2011.
\newblock \doi{10.1108/03321641111168066}.
\newblock URL \url{https://doi.org/10.1108/03321641111168066}.

\end{thebibliography}
\bibliographystyle{plainnat}
\end{document}